\documentclass[12pt]{article}
\usepackage{a4,amsfonts,amsmath,amssymb,graphicx,bbm}

\usepackage{arydshln}

\newcommand{\real}{\mathbbm{R}}

\newcommand{\im}{{\rm im\;}}
\renewcommand{\ker}{{\rm ker\;}}

\begin{document}

\newtheorem{theorem}{Theorem}
\newtheorem{lemma}{Lemma}
\newtheorem{definition}{Definition}
\newtheorem{corollary}{Corollary}
\newtheorem{condition}{Condition}


\begin{center}

{\bf \LARGE
Index-analysis for a method of lines \\[0.1ex] 
discretising multirate partial \\[0.7ex] 
differential algebraic equations}

\vspace{0.5cm}   

{\large Roland~Pulch$\mbox{}^1$\footnote{corresponding author}, Diana Est{\'e}vez Schwarz$\mbox{}^2$ and Ren{\'e} Lamour$\mbox{}^3$}

\vspace{0.1cm}

\begin{small}

{$\mbox{}^1$Institut f\"ur Mathematik und Informatik, \\
Ernst-Moritz-Arndt-Universit\"at Greifswald, \\
Walther-Rathenau-Str.~47, 
D-17489 Greifswald, Germany.}\\
Email: {\tt pulchr@uni-greifswald.de}

{$\mbox{}^2$Fachbereich II Mathematik - Physik - Chemie, \\ 
  Beuth Hochschule f\"ur Technik Berlin, \\
  Luxemburger Str.~10,
  D-13353 Berlin, Germany.}\\
Email: {\tt estevez@beuth-hochschule.de}

{$\mbox{}^3$Institut f\"ur Mathematik, 
  Humboldt-Universit\"at zu Berlin, \\
  Rudower~Chaussee~25,
  D-12489 Berlin, Germany.}\\
Email: {\tt lamour@math.hu-berlin.de}

\end{small}

\vspace{0.5cm}

\today

\end{center}



\begin{center}
{\bf Abstract}

\begin{tabular}{p{13cm}}
  In radio frequency applications, electric circuits generate signals,
  which are amplitude modulated and/or frequency modulated.
  A mathematical modelling yields typically systems of differential
  algebraic equations (DAEs).
  A multivariate signal model transforms the DAEs into multirate partial
  differential algebraic equations (MPDAEs).
  In the case of frequency modulation, an additional condition is required
  to identify an appropriate solution.
  We consider a necessary condition for an optimal solution and
  a phase condition.
  A method of lines, which discretises the MPDAEs as well as the additional
  condition, generates a larger system of DAEs. 
  We analyse the differential index of this approximative DAE system,
  where the original DAEs are assumed to be semi-explicit systems.
  The index depends on the inclusion of either differential variables or
  algebraic variables in the additional condition.
  We present results of numerical simulations for an illustrative example,
  where the index is also verified by a numerical method.
\end{tabular}
\end{center}

\clearpage


\section{Introduction}

The mathematical modelling of electric circuits often yields
time-dependent systems of nonlinear differential algebraic equations (DAEs),
see~\cite{guenther,ho,kamprent}.
In radio frequency applications, high-frequency oscillations appear,
whose amplitude and/or frequency change slowly in time.
Hence a transient simulation of initial value problems of the DAEs
is inefficient, because a numerical integrator has to capture each
oscillation.

A multidimensional signal representation yields an alternative approach.
Brachtendorf et al.~\cite{brachtbunse} derived an efficient model 
consisting of multirate partial differential algebraic equations (MPDAEs). 
Analysis and simulation of the MPDAE model in the case of
amplitude modulation without frequency modulation is given
in~\cite{bartelknorrpulch,brachtendorfetal,oliveirapedro,pulch02,roy}.

In the case of frequency modulation,
Narayan and Roychowdhury~\cite{narayan} formulated a system of
(warped) MPDAEs, where a local frequency function represents a
degree of freedom.
Hence an additional condition is required to determine a solution,
which allows for an efficient numerical simulation.
On the one hand, the MPDAEs together with phase conditions or 
similar constraints were considered 
in~\cite{pulch05,pulch06,pulch07b,pulch08proc,zhu}.
On the other hand, the identification of optimal solutions,
which exhibit a minimum amount of oscillations in some sense,
implies necessary conditions. 
MPDAEs with optimal multidimensional representations were investigated 
in~\cite{bittner-brachtendorf-a,bittner-brachtendorf-b,grebpulch,
houben-diss,houben,kugelmann-pulch,pulch08a,pulch08}.
A survey on all the above cases can be found in~\cite{pulch07a}.

In this paper, we consider initial-boundary value problems of
(warped) MPDAEs in the case of frequency modulation.
A method of lines yields a system of DAEs, whose numerical solution
approximates the exact solution.
The local frequency function is included in this approximative system.
The analytical as well as numerical properties of a DAE system are
characterised by its index, where different concepts exist,
see~\cite{griepentrog-maerz,hairer2,kunkel-mehrmann}.
We analyse the differential index of the system from the method of lines.
Therein, we assume that the circuit model consists of semi-explicit DAEs
of index one.
The focus is on the method of lines including a necessary condition
for an optimal solution from~\cite{kugelmann-pulch,pulch08a}.
In addition, the application of a phase condition going back to~\cite{narayan}
is analysed.
On the one hand, we perform a structural analysis of the DAE systems.
On the other hand, we identify the index under certain assumptions.
It follows that the index increases in most of the cases
depending on the inclusion of differential variables or algebraic
variables in the additional condition.
Furthermore, we discuss the determination of the index by
a numerical method, see~\cite{schwarz-lamour}, to confirm
our analysis.

\clearpage

The paper is organised as follows.
We review the multidimensional signal model, the MPDAE system 
and the method of lines in Section~\ref{sec:model}. 
We analyse the structure of the resulting systems of DAEs and determine
their index in Section~\ref{sec:structural} and Section~\ref{sec:index},
respectively.
Finally, Section~\ref{sec:example} depicts results from 
numerical simulations of a ring oscillator.


\section{Multirate Model}
\label{sec:model}
We review the modelling and simulation by MPDAEs.

\subsection{Problem definition}
Let the mathematical model of an electric circuit be a system of DAEs
in the general form 
\begin{equation} \label{dae}
\frac{\mbox{d}{q}({x}(t))}{\mbox{d}t} =
{f}({b}(t),{x}(t)) .
\end{equation}
The solution ${x}: \real \rightarrow \real^n$ includes unknown
voltages and currents.
The function ${b}: \real \rightarrow \real^{n_{\rm in}}$
introduces predetermined input signals.
The nonlinear functions 
${q}: \real^n \rightarrow \real^n$ and 
${f}: \real^{n_{\rm in}} \times \real^n \rightarrow \real^n$
depend on the solution.
Let ${x},{q},{f}$ be sufficiently smooth and ${b}$ be continuous.
Without loss of generality, we choose the initial time $t_0=0$.
An initial value problem is given by
\begin{equation} \label{ivp}
  {x}(0) = {x}_0
\end{equation}
with consistent initial values ${x}_0 \in \real^n$.

In radio frequency applications, the solution or some of its components
represent high-frequency oscillations. 
We assume that the input signals change slowly in a
total time interval $[0,t_{\rm end}]$.
The input signals control the amplitude and/or frequency of the solution.
Thus the function ${x}$ exhibits a huge number of oscillations
in the total time interval.
Concerning an initial value problem~(\ref{dae}), (\ref{ivp}),
a numerical integration method has to capture each oscillation by
several time steps.
Hence a transient simulation becomes inefficient due to a huge
computational effort.

\subsection{Multirate partial differential algebraic equations}
A multivariate signal model is able to decouple the slow time scale
and the fast time scale in the problem.
The solution~${x}$ of~(\ref{dae}) is represented by a multivariate function
${ \hat{x}} : [0,t_{\rm end}] \times [0,1] \rightarrow \real^n$.
The second time scale is standardised to the unit interval $[0,1]$.
The input signals~${b}$ do not require a multivariate modelling, because
they are assumed to be slowly varying functions. 
The system of DAEs~(\ref{dae}) changes into the MPDAEs, 
see~\cite[Eq.(16)]{narayan},
\begin{equation} \label{mpdae}
\frac{\partial {q}({ \hat{x}}(t_1,t_2))}{\partial t_1} + \nu(t_1)
\frac{\partial {q}({ \hat{x}}(t_1,t_2))}{\partial t_2} =
{f}({b}(t_1),{ \hat{x}}(t_1,t_2)) .
\end{equation}
Therein, the local frequency function $\nu : [0,t_{\rm end}] \rightarrow \real$
represents a degree of freedom in the multivariate modelling.
We assume that ${ \hat{x}}$ is smooth and $\nu$ is continuous.
The system~(\ref{mpdae}) is also called 'warped MPDAEs' due to 
the introduction of the local frequency function, 
which deforms the second time scale. 
The equations~(\ref{mpdae}) reveal a hyperbolic structure with a 
specific form of characteristic curves, see~\cite[Sect.4]{pulch05}.

Either initial-boundary value problems or biperiodic boundary value problems
are considered for the system of MPDAEs~(\ref{mpdae}).
In this paper, we examine initial-boundary value problems, i.e., 
\begin{equation} \label{ibvp}
{ \hat{x}}(0,t_2) = { \hat{x}}_0(t_2) , \quad 
{ \hat{x}}(t_1,t_2+1) = { \hat{x}}(t_1,t_2)
\quad \mbox{for all} \;\; t_1 \ge 0 \;\; \mbox{and} \;\; t_2 \in \real
\end{equation}
with a predetermined periodic function 
${ \hat{x}}_0 : [0,1] \rightarrow \real^n$.
The initial condition has to contain the initial values~(\ref{ivp}) by
${ \hat{x}}_0(0) = {x}_0$.
If the solution of the initial-boundary value problem has a relatively
simple form with a low amount of oscillations in the domain of definition
$[0,t_{\rm end}] \times [0,1]$,
then a numerical solution can be done efficiently.
The reason is that a coarse grid captures the multivariate function
sufficiently accurate.

A solution of the initial-boundary value problem~(\ref{mpdae}), (\ref{ibvp})
implies a solution of the original
initial value problem~(\ref{dae}), (\ref{ivp}) by
\begin{equation} \label{reconstruction}
{x}(t) := { \hat{x}} 
\left( t, \Psi(t) \right) 
\qquad \mbox{with} \qquad 
\Psi(t) := \int_0^{t} \nu(\tau) \; \mbox{d}\tau ,
\end{equation}
see~\cite[p.~902]{narayan}.
More details on the modelling by (warped) MPDAEs can be found
in~\cite{pulch07a,pulch08proc}.

\subsection{Optimal solutions}
\label{sec:optimalsolution}
In the system of MPDAEs~(\ref{mpdae}), the local frequency function
represents a degree of freedom.
The aim is to obtain a solution ${ \hat{x}}$ with a low amount of
oscillations.
Real-valued weights $w_1,\ldots,w_n \ge 0$ are considered for an
optimisation in each component.
At least one weight must be positive.
Let ${ W} := {\rm diag}(w_1,\ldots,w_n) \in \real^{n \times n}$.
In~\cite{kugelmann-pulch}, the minimisation of the functional
\begin{equation} \label{functional}
  J({ \hat{x}}) = \int_0^{T_1} \int_0^{T_2}
  \left\| { W}^{\frac{1}{2}} \textstyle
  \frac{\partial { \hat{x}}}{\partial t_2}
  \right\|_2^2 \; {\rm d}t_2 {\rm d}t_1
\end{equation}
was examined with
${ W}^{\frac{1}{2}} := {\rm diag} \left( \sqrt{w_1},\ldots,\sqrt{w_n} \right)$
and the Euclidean norm $\| \cdot \|_2$.
This optimisation turns out to be equivalent to the point-wise
minimisation of the functional
\begin{equation} \label{functional-pw}
  \tilde{J}( t_1 ; { \hat{x}} ) =
  \int_0^{T_2}
  \left\| { W}^{\frac{1}{2}} \textstyle
  \frac{\partial { \hat{x}}}{\partial t_2}
  \right\|_2^2 \; {\rm d}t_2
  \qquad \mbox{for each} \;\; t_1 .
\end{equation}
In our case, it holds that $T_1 = t_{\rm end}$ and $T_2=1$.
Existence and uniqueness of optimal solutions was also proven
in~\cite{kugelmann-pulch}.

A necessary condition for a solution, which minimises both the
functional~(\ref{functional}) and~(\ref{functional-pw}),
reads as
\begin{equation} \label{opt-necessary}
  \int_0^1
  \left( \textstyle \frac{\partial { \hat{x}}}{\partial t_1} \right)^\top
  { W} \left( \textstyle \frac{\partial { \hat{x}}}{\partial t_2} \right)
  \; \mbox{d}t_2 = 0
  \qquad \mbox{for all}\;\; t_1 \in [0,t_{\rm end}] ,
\end{equation}
see~\cite[Cor.~1]{kugelmann-pulch}.
This constraint was already found for an equivalent optimisation criterion
in~\cite{pulch08a}. 

\subsection{Phase conditions}
\label{sec:phasecondition}
Alternatively, phase conditions can be added to the
MPDAEs~(\ref{mpdae})
either in the time domain or in the frequency domain,
see~\cite{narayan,zhu}.
In the time domain, continuous phase conditions just represent
an additional boundary condition at $t_2=0$.
A component $\ell \in \{ 1,\ldots,n \}$ has to be chosen from
${ \hat{x}} = (\hat{x}_1,\ldots,\hat{x}_n)^\top$.
In~\cite[Eq.(9)]{narayan}, a derivative of this component is
predetermined as a slowly varying function~$\eta$.
This approach can be written as
$$ \textstyle \frac{\partial \hat{x}_{\ell}}{\partial t_2} (t_1,0) =
\eta(t_1) \qquad \mbox{for all}\;\; t_1 \in [0,t_{\rm end}] . $$
In~\cite[Eq.(7)]{pulch06}, the function~$\hat{x}_{\ell}$ itself
is forced to be a constant value at the boundary.
Hence we investigate the phase condition
\begin{equation} \label{phasecondition}
  \hat{x}_{\ell}(t_1,0) = \eta(t_1)
  \qquad \mbox{for all}\;\; t_1 \in [0,t_{\rm end}] 
\end{equation}
with a predetermined function $\eta : [0,t_{\rm end}] \rightarrow \real$.
As mentioned above,
often a constant choice $\eta \equiv \eta_0$ is feasible.

There is a heuristic motivation of the phase conditions.
If a component of the solution exhibits a simple slowly varying
shape on the boundary, then most likely the complete solution has
an elementary behaviour with a low amount of oscillations.
In comparison to a minimisation of the functional~(\ref{functional}),
phase conditions often yield suboptimal solutions.

\subsection{Method of lines}
In a method of lines, the second derivative of the MPDAE system~(\ref{mpdae})
is replaced by finite differences.
This discretisation is applied on the lines $t_{2,i} := (i-1)h$
for $i = 1,\ldots,m$ with a step size $h = \frac{1}{m}$
given some integer~$m$.
We obtain a system of DAEs
\begin{equation} \label{mol}
  \frac{{\rm d} {q}({ \bar{x}}_i(t_1))}{{\rm d} t_1} =
  {f}({b}(t_1),{ \bar{x}}_i(t_1)) - \nu(t_1)
  ({\rm D}_i ({q},{ \bar{x}})) (t_1) 
\end{equation}
for $i=1,\ldots,m$.
The numerical solution is
${ \bar{x}} := ({ \bar{x}}_1^\top , \ldots , { \bar{x}}_m^\top)^\top
\in \real^{mn}$.
Each function ${ \bar{x}}_i (t_1)$ represents an approximation of
the solution ${ \hat{x}}(t_1,t_{2,i})$ for $i=1,\ldots,m$.
The symbol ${\rm D}_i$ denotes the finite difference formula.
For example, the backward differentiation formulas (BDF)
of order one and two, see~\cite[Ch.III.1]{hairer1},
read as
\begin{equation} \label{bdf}
  \begin{array}{rcl}
    {\rm D}_i ({q},{ \bar{x}}) & = &
    \frac{1}{h}
    \left[ {q}({ \bar{x}}_i)
      - {q}({ \bar{x}}_{i-1}) \right] , \\[1ex]
    {\rm D}_i ({q},{ \bar{x}}) & = &  
    \frac{1}{h}
    \left[ \frac{3}{2} {q}({ \bar{x}}_i) - 2 {q}({ \bar{x}}_{i-1})
      + \frac{1}{2} {q}({ \bar{x}}_{i-2}) \right]  \\
  \end{array}
\end{equation}
for $i=1,\ldots,m$.
Therein, the periodicities ${ \bar{x}}_j = { \bar{x}}_{j+m}$
for each~$j$ have to be used to eliminate the unknowns
for $j \notin \{ 1,\ldots,m \}$.
The usual choices of finite difference approximations are convergent
for sufficiently smooth functions.
The system~(\ref{mol}) is still underdetermined,
because an appropriate local frequency function~$\nu$ is not
identified yet.
We require an additional condition.

Now our aim is to determine an optimal solution as introduced
in Section~\ref{sec:optimalsolution}.
In a numerical method, the necessary condition~(\ref{opt-necessary})
has to be discretised using the lines.
Firstly, the integral is replaced by the rectangular rule,
which is equivalent to the trapezoidal rule due to the periodicity in~$t_2$.
Furthermore, the derivative with respect to~$t_2$ is approximated by a
finite difference scheme.
We obtain
\begin{equation} \label{opt-disc}
  \sum_{i=1}^m 
  \left( \textstyle \frac{{\rm d} { \bar{x}}_i}{{\rm d} t_1} \right)^\top
  { W} \left( {\rm D}_i ({ \bar{x}}) \right) = 0
  \qquad \mbox{for all}\;\; t_1 \in [0,t_{\rm end}] ,
\end{equation}
where we apply the same finite differences as in~(\ref{mol}).
More precisely, it holds that
${\rm D}_i ({ \bar{x}}) = {\rm D}_i ({q},{ \bar{x}})$
using the identity function as~${q}$.
The system~(\ref{mol}), (\ref{opt-disc}) involves as many equations
as unknowns.
For a linear function ${q}({x}) = Q {x}$ with a
constant mass matrix $Q \in \real^{n \times n}$, the complete system
can be written in the quasi-linear form
\begin{equation} \label{dae-quasi-linear}
  { M}({ y}) \dot{ y} = {f}(t_1,{ y})
\end{equation}
with ${ y} = ({ \bar{x}}^\top,\nu)^\top$ and
a state-dependent mass matrix
\begin{equation} \label{massmatrix1}
{ M}({ y}) =
\left(
\begin{array}{ccc:c}
  Q & & & { 0} \\
  & \ddots & & \vdots \\
  & & Q & { 0} \\ \hdashline
  ({ W} {\rm D}_1({ \bar{x}}))^\top & \cdots &
  ({ W} {\rm D}_m({ \bar{x}}))^\top \rule{0ex}{2.5ex}
  & 0 \\
\end{array}
\right) \in \real^{(mn+1) \times (mn+1)}.
\end{equation}
Alternatively, a phase condition from Section~\ref{sec:phasecondition}
can be included in the method of lines.
The boundary $t_2=0$ coincides with the line $t_{2,i}$ for $i=1$.
Choosing a component $\ell \in \{ 1,\ldots,n \}$,
the phase condition~(\ref{phasecondition}) yields
\begin{equation} \label{phase-mol}
  \bar{x}_{1,\ell} (t_1) = \eta(t_1)
  \qquad \mbox{for all}\;\; t_1 \in [0,t_{\rm end}] 
\end{equation}
with a predetermined function~$\eta$.
The system~(\ref{mol}), (\ref{phase-mol}) has as many equations
as unknowns again.
In the case of a linear function~${q}$, the system features
the form~(\ref{dae-quasi-linear}) with the constant mass matrix
\begin{equation} \label{massmatrix2}
{ M} =
\left(
\begin{array}{ccc:c}
  Q & & & { 0} \\
  & \ddots & & \vdots \\
  & & Q & { 0} \\ \hdashline
  0 & \cdots & 0 & 0 \\
\end{array}
\right) \in \real^{(mn+1) \times (mn+1)} .
\end{equation}
The mass matrices~(\ref{massmatrix1}), (\ref{massmatrix2}) are always
singular due to the last column.
It follows that also a system~(\ref{dae}) of ordinary differential equations
changes to DAEs in the method of lines.

\subsection{Semi-explicit systems}
In the following sections,
we restrict the analysis to the important case of semi-explicit DAEs.
Such systems are characterised by a linear function $q(x) = Dx$
with a diagonal matrix
${ D} = {\rm diag}(1,\ldots,1,0,\ldots,0)$.
The solution ${x} = ({ y}^\top,{ z}^\top)^\top$
is partitioned into differential variables ${ y} \in \real^{n_y}$
and algebraic variables ${ z} \in \real^{n_z}$ ($n_y+n_z=n$).
Now the system~(\ref{dae}) reads as 
\begin{equation} \label{semi-expl-dae}
  \begin{array}{rcl}
    \dot{ y}(t) & = & {f}(t,{ y}(t),{ z}(t)) , \\[1ex]
    { 0} & = & { g}(t,{ y}(t),{ z}(t)) . \\
  \end{array}
\end{equation}
The dependence on the input signals ${b}$ is represented by
the first argument of the functions ${f}$ and ${ g}$ for
notational convenience.
  
If a system~(\ref{dae}) is given including a linear function
${q}({x}) = Q {x}$
with a constant mass matrix~$Q \in \real^{n \times n}$,
then it can be transformed into a semi-explicit system of DAEs.
For example, a singular value decomposition of the matrix~$Q$
can be used for this transformation.

The exact definition of the differential index can be found
in~\cite[Ch.VI.5]{hairer2}, for example.
Roughly speaking, the differential index is the minimum number of
differentiations applied to a DAE system such that an
ordinary differential equation (ODE) can be derived for all
variables of the solution.

The mathematical modelling of electric circuits yields typically
systems of DAEs with differential index either one or two.
We restrict the analysis to semi-explicit systems~(\ref{semi-expl-dae})
of index one.
In this case, the Jacobian matrix
$\frac{\partial { g}}{\partial { z}}$
is always non-singular.
Consequently, both the differential index and the perturbation index
of~(\ref{semi-expl-dae}) are equal to one.
It follows that we achieve an ODE for the algebraic variables after
one differentiation of the system.

We consider the transition from the DAEs~(\ref{semi-expl-dae})
to the MPDAEs
\begin{equation} \label{mpdae-semi-expl}
 \begin{array}{rcl} 
   \frac{\partial  \hat{y}}{\partial t_1} + \nu(t_1)
   \frac{\partial  \hat{y}}{\partial t_2}
   & = & {f}(t_1,{ \hat{y}},{ \hat{z}}) , \\[1ex]
    { 0} & = & { g}(t_1,{ \hat{y}},{ \hat{z}}) . \\
  \end{array}
\end{equation}
Now the initial-boundary value problem~(\ref{ibvp}) reads as
\begin{equation} \label{ibvp-semi-expl}
  \begin{array}{rclrcl}
    { \hat{y}}(t_1,0) & = & { \hat{y}}_0(t_1) , \quad &
    { \hat{y}}(t_1,t_2) & = & { \hat{y}}(t_1,t_2+1) , \\[0.5ex]
    { \hat{z}}(t_1,0) & = & { \hat{z}}_0(t_1) , \quad & 
    { \hat{z}}(t_1,t_2) & = & { \hat{z}}(t_1,t_2+1) \\
  \end{array}
\end{equation}
for all~$t_1 \ge 0$ and all $t_2 \in \real$.
The system~(\ref{mol}) from the method of lines becomes
\begin{equation} \label{mol-semi-expl}
  \begin{array}{rcl}
    {\dot{ { \bar{ y}} }}_i (t) & = &
        {f}(t,{ \bar{y}}_i(t),{ \bar{z}}_i(t)) -
        \nu(t) ( {\rm D}_i ({ \bar{y}}) ) (t) , \\[1ex]
    { 0} & = & { g}(t,{ \bar{y}}_i(t),{ \bar{z}}_i(t))  \\
  \end{array}
\end{equation}
for $i=1,\ldots,m$.
For notational convenience, the slow time variable~$t_1$ is replaced by~$t$
and differentiations are indicated by a dot.
A general finite difference formula exhibits the structure
\begin{equation} \label{difference-formula}
  {\rm D}_i ({ \bar{y}}) =
  \frac{1}{h} \sum_{j=-q}^p \alpha_j { \bar{y}}_{i+j}
\end{equation}
for $i=1,\ldots,m$ with real coefficients $\alpha_{-q},\ldots,\alpha_p$
and integers $q,p \ge 0$.
Thus a differentiation just reads as
$$ {\textstyle \frac{{\rm d}}{{\rm d}t}}
 {\rm D}_i ({ \bar{y}}) =
 \frac{1}{h} \sum_{j=-q}^p \alpha_j {\dot{ { \bar{ y}} }}_{i+j} =
 \frac{1}{h} \sum_{j=-q}^p \alpha_j \left[
 {f}(t,{ \bar{y}}_{i+j},{ \bar{z}}_{i+j}) -
        \nu(t) ( {\rm D}_{i+j} ({ \bar{y}}) ) \right] $$
 for $i=1,\ldots,m$.

In the case of semi-explicit systems~(\ref{mol-semi-expl}),
the phase condition~(\ref{phase-mol}) becomes either
\begin{equation} \label{phase-semi-expl}
  \bar{y}_{1,\ell}(t) = \eta(t)
  \qquad \mbox{or} \qquad
  \bar{z}_{1,\ell}(t) = \eta(t)
  \qquad \mbox{for all} \;\; t
\end{equation}
with either an $\ell \in \{ 1,\ldots,n_y \}$ or
an $\ell \in \{ 1, \ldots , n_z \}$, respectively.

Considering the functional~(\ref{functional}) of the optimisation, 
the weights are partitioned now into $w_1^y,\ldots,w_{n_y}^y$ and
$w_1^z,\ldots,w_{n_z}^z$ for differential variables and algebraic variables,
respectively.
The necessary condition~(\ref{opt-disc}) results to
\begin{equation} \label{opt-disc2}
\displaystyle \sum_{i=1}^m \bigg( \sum_{\ell = 1}^{n_y} 
w_\ell^y \cdot \dot{\bar{y}}_{i,\ell} \cdot {\rm D}_{i,\ell} ({ \bar{y}})
+ \sum_{\ell = 1}^{n_z} 
w_\ell^z \cdot \dot{\bar{z}}_{i,\ell} \cdot {\rm D}_{i,\ell} ({ \bar{z}})
\bigg) = 0 
\end{equation}
for all $t \in [0,t_{\rm end}]$.


\newpage
 
\section{Structural Analysis}
\label{sec:structural}
We examine the general structure of the semi-explicit systems of DAEs,
which are generated by method of lines and additional conditions.

\subsection{General setting}
Concerning the semi-explicit system~(\ref{mol-semi-expl}),
we introduce the vectors
\begin{eqnarray*}
x_1&:=&(\bar{y}_{1,1}, \ldots,\bar{y}_{1,n_y}, \ldots, \bar{y}_{m ,1}, \ldots, \bar{y}_{m,n_y})^\top , \\
x_2&:=&(\bar{z}_{1,1}, \ldots, \bar{z}_{1,n_z}, \ldots, \bar{z}_{m ,1}, \ldots,\bar{z}_{m,n_z})^\top , \\
x_3&:=&\nu , 
\end{eqnarray*}
where
$x_1 : \real \rightarrow \real^{n_1}$ with $n_1 := m n_y$ and
$x_2 : \real \rightarrow \real^{n_2}$ with $n_2 := m n_z$
include the differential variables and the algebraic variables, respectively,
and $x_3$ denotes the local frequency function.
Let $\bar{n} := n_1 + n_2 + 1$.
In the method of lines, each system of DAEs exhibits the general form
$f(t,x,\dot{x})=0$
with the variables $x=(x_1^\top,x_2^\top,x_3)^\top$.
More detailed, we obtain the structure
\begin{eqnarray}
\dot{x}_1 & = & f_1(t,x_1,x_2,x_3),  \label{eq:1DAEstructure}\\
0 & = & f_2(t,x_1,x_2),\label{eq:2DAEstructure}\\
0 & = & f_3(t,x_1,x_2,\dot{x}_1,\dot{x}_2),\label{eq:3DAEstructure}
\end{eqnarray}
with $f_1 \in \real^{n_1}$, $f_2 \in \real^{n_2}$ from the semi-discretisation
and $f_3 \in \real$ from an additional condition.

For our consideration of the derivative array \cite{BrCaPe96}, we define
\begin{eqnarray} \label{eq:DAE_for_x_deriv}
  \hat{f}(t,x,\dot{x},\ddot{x}):=
  \textstyle \frac{{\rm d}}{{\rm d}t}f(t,x,\dot{x})
\end{eqnarray}
and likewise $\hat{f}_1$, $\hat{f}_2$, $\hat{f}_3$.
It holds that
$\frac{\partial f_i}{\partial x_j} = \frac{\partial \hat{f}_i}{\partial \dot{x}_j}$
for $i,j=1, 2$ due to the chain rule of differentiation. 

In the following, we analyse the structure and the index of the
DAE system \eqref{eq:1DAEstructure}-\eqref{eq:3DAEstructure} as
described in \cite{schwarz-lamour2016} and \cite{schwarz-lamour}.
Thus matrices
${\cal B}^{[k]} \in \real^{\bar{n}(k+1) \times \bar{n}(k+1)}$
involving
$\frac{\partial f}{\partial x}$, $\frac{\partial f}{\partial \dot{x}}$,  $\frac{\partial \hat{f}}{\partial x}$, $\frac{\partial \hat{f}}{\partial \dot{x}}$, $\frac{\partial \hat{f}}{\partial \ddot{x}}$
are considered.
In fact, we will check if the matrices ${\cal B}^{[k]}$ are 1-full
with respect to the first $\bar{n}$ columns for $k=1,2$,
i.e., whether
\begin{equation} \label{kernel-B}
\ker {\cal B}^{[k]} \subseteq \left\{ \begin{pmatrix}
	s_0 \\
	s_1 \\
\end{pmatrix}
\; : \; s_0 \in \real^{\bar{n}} , \; s_0 =0 , \; s_1 \in \real^{\bar{n}k} \right\}.
\end{equation}
As a result, we obtain index criteria
and a characterisation of the higher-index component. 
In fact, the 1-fullness we check characterises if we can represent $x_3=\nu$ 
as a function of $(x_1,x_2)$ directly ($k=1$ and index one) or after one differentiation ($k=2$ and index two).
A differentiation of this function would deliver the expression for $\dot{x}_3=\dot{\nu}$.
In \cite{schwarz-lamour2016} it has been shown that in case that this index is defined, it
coincides with the differential index. Roughly speaking, this means that the
relevant subspaces related to the DAE structure have constant dimensions, in 
accordance to the concepts of \cite{LamourMaerzTischendorf2013}.
Otherwise, there may be singular points.
In \cite{schwarz-lamour2015a,schwarz-lamour2015b},
the equivalence of this index based on 1-fullness and the tractability
index from \cite{LamourMaerzTischendorf2013} has been proofed for wide classes of DAEs of
index up to two.

We will repeatedly make use of the fact that for bordered matrices of the form
\[
\left(
\begin{array}{c:c}
	v & H \\ \hdashline 
	0 & w \\
\end{array} 
\right), \quad H \in \real^{(\bar{n}-1) \times (\bar{n}-1)} \; \mbox{non-singular}, \quad v \in \real^{(\bar{n}-1) \times 1}, \quad w \in \real^{1 \times (\bar{n}-1)} ,
\]
it holds that the bordered matrix is 1-full with respect to the first column,
i.e.,
\[
\ker \begin{pmatrix}
	v & H \\
	0 & w
\end{pmatrix} \subseteq \left\{ \begin{pmatrix}
	s_0 \\
	s_1
\end{pmatrix}
\; : \; s_0 \in \real, \; s_0=0 , \;  s_1 \in \real^{\bar{n}-1} \right\} ,
\]
if and only if the bordered matrix is non-singular, i.e., iff $w H^{-1}v \neq 0$.

In all our considerations, we assume that
$
\frac{\partial f_2}{\partial x_2} = \frac{\partial \hat{f}_2}{\partial x'_2}
$
is non-singular and $\frac{\partial f_1}{\partial x_3} \neq 0$. 
The structure of the DAE will particularly depend on the additional
equation \eqref{eq:3DAEstructure}, i.e.,
either \eqref{phase-semi-expl} or \eqref{opt-disc2},
because it determines the vector $w$ of the bordered matrices. 
In contrast, the vector $v$ depends on $\frac{\partial f_1}{\partial x_3}$.

\subsection{DAEs for phase condition}
\label{sec:structure-phase}
We consider \eqref{mol-semi-expl} together with \eqref{phase-semi-expl},
where the system \eqref{eq:1DAEstructure}-\eqref{eq:3DAEstructure}
presents the structure
\begin{eqnarray*}
\dot{x}_1&=&f_{1,1}(t,x_1,x_2)+F_{1,3}(x_1)x_3, \label{eq:1DAEstructurePhaseCond}\\
0 &=& f_2(t,x_1,x_2),\label{eq:2DAEstructurePhaseCond}\\
0 &=& F_{3,1}x_1+F_{3,2}x_2 +f_3(t)\label{eq:3DAEstructurePhaseCond}
\end{eqnarray*}
with constant matrices $F_{3,1} \in \real^{1\times n_1}$, $F_{3,2} \in \real^{1\times n_2}$ satisfying either
\[
F_{3,1} \neq \left(0, \ldots, 0\right) \quad \mbox{and} \quad F_{3,2}=\left(0, \ldots, 0\right)
\]
or
\[
F_{3,1}= \left(0, \ldots, 0\right) \quad \mbox{and} \quad F_{3,2}\neq \left(0, \ldots, 0\right).
\]
For $A:= \frac{\partial f}{\partial \dot{x}}$ we obtain 
\[
A=\begin{pmatrix}
	I & 0 & 0 \\
	0 & 0 & 0 \\
	0 & 0 & 0 
\end{pmatrix}, \quad Q = \begin{pmatrix}
	0 & 0 & 0 \\
	0 & I & 0 \\
	0 & 0 & 1
\end{pmatrix}, \quad P=A ,
\]
where $Q$ and $P=I-Q$ denote the orthogonal projectors onto $\ker A$
and $(\ker A)^\perp$, respectively.

\subsubsection{Check of index-1 condition}
\label{sec:structure-phase-index1}
 In the first step, we consider the matrix
\begin{eqnarray}
{\cal B}^{[1]}=\begin{pmatrix}
	I & 0 & 0  & & &\\
	0 & 0 & 0  & & &\\
	0 & 0 & 0  & & &\\
	\frac{\partial f_1}{\partial x_1} & \frac{\partial f_1}{\partial x_2} & F_{1,3} & I & 0 & 0\\[0.5ex]
	\frac{\partial f_2}{\partial x_1} & \frac{\partial f_2}{\partial x_2}& 0&0 &0 & 0 &\\[0.5ex]
	F_{3,1}                           & F_{3,2}& 0 & 0 & 0 & 0 
\end{pmatrix}.
\label{eq:calB1phasecond}
\end{eqnarray}
In terms of the index-definition from \cite{schwarz-lamour}, the index is one if and only if ${\cal B}^{[1]}$ 
is 1-full with respect to the first $\bar{n}$ columns, i.e., if
\[
\ker {\cal B}^{[1]} \subseteq \left\{ \begin{pmatrix}
	s_0 \\
	s_1
\end{pmatrix} \in \real^{2\bar{n}} \; : \; s_0 =0 \right\} , 
\]
cf.~(\ref{kernel-B}).
Note that the full rank of $\frac{\partial f_2}{\partial x_2}$ implies
\[
\ker {\cal B}^{[1]}= \ker \begin{pmatrix}
	I & 0 & 0  & & &\\
	0 & I & 0  & & &\\
	0 & 0 & F_{1,3} & I & 0 & 0
	\end{pmatrix}.
\]
Therefore, ${\cal B}^{[1]}$ is not 1-full and, consequently,
the index cannot be one.
We also recognise that $x_1$ and $x_2$ are not higher-index variables.
The orthogonal projector $T \in \real^{\bar{n} \times \bar{n}}$,
which describes the higher-index component $x_3$,
reads as
\begin{equation} \label{projector-T}
T=\begin{pmatrix}
	0 & 0 & 0\\
	0 & 0 & 0 \\
	0 & 0 & 1
\end{pmatrix}.
\end{equation}

\subsubsection{Check of index-2 condition}
In the case that the index of the DAE was not one,
according to \cite{schwarz-lamour}, it is two 
if and only if the matrix
\[
{\cal B}^{[2]}=\begin{pmatrix}
	I & 0 & 0  & & & & & &\\
	0 & 0 & 0  & & & & & &\\
	0 & 0 & 0  & & & & & &\\[0.5ex]
	\frac{\partial f_1}{\partial x_1} & \frac{\partial f_1}{\partial x_2} & F_{1,3} & I & 0 & 0\\[0.5ex]
	\frac{\partial f_2}{\partial x_1} & \frac{\partial f_2}{\partial x_2}& 0&0 &0 & 0 &\\[0.5ex]
F_{3,1}&F_{3,2} & 0 & 0 & 0 & 0\\[0.5ex]
  \frac{\partial \hat{f}_1}{\partial x_1} & \frac{\partial \hat{f}_1}{\partial x_2} &  0& \frac{\partial \hat{f}_1}{\partial \dot{x}_1} & \frac{\partial \hat{f}_1}{\partial \dot{x}_2} & F_{1,3} & I & 0 & 0\\[0.5ex]
	\frac{\partial \hat{f}_2}{\partial x_1} & \frac{\partial \hat{f}_2}{\partial x_2} &  0&\frac{\partial \hat{f}_2}{\partial \dot{x}_1} & \frac{\partial \hat{f}_2}{\partial \dot{x}_2}& 0& 0 &0 & 0 &\\[0.5ex]
	0 & 0 & 0 & F_{3,1}&F_{3,2} & 0 &0 &0 &0
\end{pmatrix}
\]
is 1-full with respect to the first $\bar{n}$ columns.
If this property is not satisfied, then the index may be higher than two
or not defined.

Under the assumption that $\frac{\partial f_2}{\partial x_2}=\frac{\partial \hat{f}_2}{\partial \dot{x}_2}$ is non-singular, we now obtain
\[
\ker {\cal B}^{[2]} = \ker
\begin{pmatrix}
	I & 0 & 0  & & & & & &\\
	0 & I & 0  & & & & & &\\
%
	0 & 0 & F_{1,3} & I & 0 & 0\\
%
	0 & 0 &  0&\frac{\partial \hat{f}_2}{\partial \dot{x}_1} & \frac{\partial \hat{f}_2}{\partial \dot{x}_2}& 0& 0 &0 & 0 &\\
	0 & 0 & 0 & F_{3,1}&F_{3,2} & 0 &0 &0 &0\\
	0 & 0 &  0& \frac{\partial \hat{f}_1}{\partial \dot{x}_1} & \frac{\partial \hat{f}_1}{\partial \dot{x}_2} & F_{1,3} & I & 0 & 0
\end{pmatrix} 
\]
and deduce
\begin{enumerate}
	\item Assuming
	$F_{3,1} \neq  \left(0, \ldots, 0\right)$ and $F_{3,2}=\left(0, \ldots, 0\right)$, 
${\cal B}^{[2]}$ is 1-full if and only if 
\[
\begin{pmatrix}
F_{1,3} & I \\
0 & F_{3,1}	
\end{pmatrix}
\]
is non-singular. Consequently, the index is two if and only if
\begin{equation} \label{ftimesf}
  F_{3,1} \cdot F_{1,3} \neq 0 .
\end{equation}
\item
  Assuming $F_{3,1}= \left(0, \ldots, 0\right)$ and $F_{3,2}\neq \left(0, \ldots, 0\right)$, 
${\cal B}^{[2]}$ is 1-full if and only if the bordered matrix
\[
\begin{pmatrix}
	F_{1,3} & I & 0  \\
	0& \frac{\partial \hat{f}_2}{\partial \dot{x}_1} & \frac{\partial \hat{f}_2}{\partial \dot{x}_2} \\
	0 & 0 &F_{3,2} 
\end{pmatrix}=
\left(
\begin{array}{c:cc}
	F_{1,3} & I & 0  \\
	0& \frac{\partial f_2}{\partial x_1} & \frac{\partial f_2}{\partial x_2} \\[0.5ex] \hdashline 
	0 & 0 &F_{3,2} 
\end{array} 
\right)
\]
is non-singular.
Thus the index turns out to be two if and only if
\[
\begin{pmatrix}
	 \frac{\partial f_2}{\partial x_1} \cdot F_{1,3}& \frac{\partial f_2}{\partial x_2} \\
	 0 &F_{3,2} 
\end{pmatrix}
\]
is non-singular, i.e., iff $F_{3,2} \cdot  \left( \frac{\partial f_2}{\partial x_2}\right)^{-1} \cdot\frac{\partial f_2}{\partial x_1} \cdot F_{1,3} \neq 0$.

\end{enumerate}

If the above assumptions on non-singularity are not fulfilled,
then the index may be higher or not defined. 

\subsection{DAEs for optimal solutions}
\label{sec:structural-opt}
For \eqref{mol-semi-expl} together with \eqref{opt-disc2} the 
structure of the system \eqref{eq:1DAEstructure}-\eqref{eq:3DAEstructure}
reads as
\begin{eqnarray}
\dot{x}_1&=&f_{1,1}(t,x_1,x_2)+F_{1,3}(x_1)x_3, \\
0 &=& f_2(t,x_1,x_2),\\
0 &=& F_{3,1}(x_1) \dot{x}_1 + F_{3,2}(x_2) \dot{x}_2, \label{eq:opt-general}
\end{eqnarray}
with $F_{1,3} : \real^{n_1} \rightarrow \real^{n_1}$,
$F_{1,3} \neq \left(0, \ldots,0\right)^\top$,
$F_{3,1}: \real^{n_1} \rightarrow \real$, and
$F_{3,2}: \real^{n_2} \rightarrow \real$.
For $A:= \frac{\partial f}{\partial \dot{x}}$, it follows that
\[
A(x_1,x_2) =\begin{pmatrix}
	I & 0 & 0 \\
	0 & 0 & 0 \\
	F_{3,1}(x_1) & F_{3,2}(x_2) & 0 
\end{pmatrix} .
\]
Let $Q_{3,2}$ be the orthogonal projector onto $\ker F_{3,2}(x_2)$.
Hence the orthogonal projector $Q$ onto $\ker A$ results to
\[
Q(x_1,x_2)=\begin{pmatrix}
	0 & 0 & 0 \\
	0 & Q_{3,2} & 0 \\
	0 & 0 & 1
\end{pmatrix}.
\]
Again, we determine the index by rank considerations.

\subsubsection{Check of index-1 condition}

Let $P_{3,2}=I-Q_{3,2}$.
The index will be one, if and only if
\[
{\cal B}^{[1]}=\begin{pmatrix}
	I & 0 & 0  & & &\\
	0 & P_{3,2} & 0  & & &\\
	0 & 0 & 0  & & &\\[0.5ex]
	\frac{\partial f_1}{\partial x_1} & \frac{\partial f_1}{\partial x_2} & F_{1,3} & I & 0 & 0\\[0.5ex]
	\frac{\partial f_2}{\partial x_1} & \frac{\partial f_2}{\partial x_2}& 0&0 &0 & 0 &\\[0.5ex]
	\frac{\partial f_3}{\partial x_1} & \frac{\partial f_3}{\partial x_2} & 0 &F_{3,1}&F_{3,2} & 0 
\end{pmatrix}
\]
is 1-full with respect to the first $\bar{n}$ columns.
Since we assume that  $\frac{\partial f_2}{\partial x_2}$ is non-singular, it obviously holds that
\[
\ker {\cal B}^{[1]} = \ker \begin{pmatrix}
	I & & & & & \\
	& I & & & & \\
	  & & F_{1,3} & I & & \\
    & & &	F_{3,1} & F_{3,2} & 0 \\
\end{pmatrix}.
\]
Again $x_1$ and $x_2$ are not higher-index variables in this situation.

Let us now analyse different cases:
\begin{enumerate}
	\item 
If $F_{3,2}=0$, then 1-fullness is given if and only if
\[
\begin{pmatrix}
	F_{1,3} & I \\
	0 & F_{3,1}
\end{pmatrix}
\]
is non-singular.
This condition is obviously fulfilled if and only if
it holds that~(\ref{ftimesf}). 
Consequently, for $F_{3,2}=0$, the index is one if and only if
the property~(\ref{ftimesf}) is satisfied.

\item
  If $F_{3,2}\neq 0$, then ${\cal B}^{[1]}$ is not 1-full and the index is
  higher than one or may be not defined. 
  This property follows from the fact that for $F_{3,2}\neq 0$ the matrix
\[
\begin{pmatrix}
	I & 0 \\
	F_{3,1} & F_{3,2}
\end{pmatrix}
\]
has full row rank $n_1+1$ and therefore
\[
\begin{pmatrix}
F_{1,3}\\
0
\end{pmatrix}
 \in \im \begin{pmatrix}
	I & 0 \\
	F_{3,1} & F_{3,2}
\end{pmatrix} =\real^{n_1+1}.
\]
Consequently, ${\cal B}^{[1]}$ is not 1-full.
\end{enumerate}

\subsubsection{Check of index-2 condition}
If the index is one,
the next step in the index analysis is to check whether
\[
{\cal B}^{[2]}=\begin{pmatrix}
	I & 0 & 0  & & & & & &\\
	0 & P_{3,2} & 0  & & & & & &\\
	0 & 0 & 0  & & & & & &\\[0.5ex]
	\frac{\partial f_1}{\partial x_1} & \frac{\partial f_1}{\partial x_2} & F_{1,3} & I & 0 & 0\\[0.5ex]
	\frac{\partial f_2}{\partial x_1} & \frac{\partial f_2}{\partial x_2}& 0&0 &0 & 0 &\\[0.5ex]
	\frac{\partial f_3}{\partial x_1} & \frac{\partial f_3}{\partial x_2} & 0 &F_{3,1}&F_{3,2} & 0 \\[0.5ex]
  \frac{\partial \hat{f}_1}{\partial x_1} & \frac{\partial \hat{f}_1}{\partial x_2} &  0& \frac{\partial \hat{f}_1}{\partial \dot{x}_1} & \frac{\partial \hat{f}_1}{\partial \dot{x}_2} & F_{1,3} & I & 0 & 0\\[0.5ex]
	\frac{\partial \hat{f}_2}{\partial x_1} & \frac{\partial \hat{f}_2}{\partial x_2} &  0&\frac{\partial \hat{f}_2}{\partial \dot{x}_1} & \frac{\partial \hat{f}_2}{\partial \dot{x}_2}& 0& 0 &0 & 0 &\\[0.5ex]
	\frac{\partial \hat{f}_3}{\partial x_1} & \frac{\partial \hat{f}_3}{\partial x_2} & 0 &\frac{\partial \hat{f}_3}{\partial \dot{x}_1} & \frac{\partial \hat{f}_3}{\partial \dot{x}_2} & 0 &F_{3,1}&F_{3,2} & 0 
\end{pmatrix}
\]
is 1-full with respect to the first $\bar{n}$ columns.
If this criterion is not satisfied, then the index may be larger than two
or not defined.
Since $\frac{\partial f_2}{\partial x_2}$ is invertible,
it holds that
\[
\ker {\cal B}^{[2]} = \ker 
\begin{pmatrix}
	I & 0 & 0  & & & & & &\\
	0 & I & 0  & & & & & &\\
%
	0 & 0 &F_{1,3} & I & 0 & 0\\
	0 & 0 & 0 &F_{3,1}&F_{3,2} & 0 \\[0.5ex]
%
0 & 0 &  0& \frac{\partial \hat{f}_1}{\partial \dot{x}_1} & \frac{\partial \hat{f}_1}{\partial \dot{x}_2} & F_{1,3} & I & 0 & 0\\[0.5ex]
	0 & 0 &0&\frac{\partial \hat{f}_2}{\partial \dot{x}_1} & \frac{\partial \hat{f}_2}{\partial \dot{x}_2}& 0& 0 &0 & 0 &\\[0.5ex]
	0 & 0 & 0 &\frac{\partial \hat{f}_3}{\partial \dot{x}_1} & \frac{\partial \hat{f}_3}{\partial \dot{x}_2} & 0 &F_{3,1}&F_{3,2} & 0 
\end{pmatrix}.
\]

Consequently, ${\cal B}^{[2]}$ is 1-full with respect to the first $\bar{n}$ columns, if and only if
the matrix
\begin{eqnarray}
\label{eq:ReducedMatrixInd2OptimalSol}
\begin{pmatrix}
   F_{1,3} & I & 0 & 0\\[0.5ex]
   0 &F_{3,1}&F_{3,2} & 0 \\[0.5ex]
	 0 &\frac{\partial \hat{f}_2}{\partial \dot{x}_1} & \frac{\partial \hat{f}_2}{\partial \dot{x}_2}& 0& 0 &0 & \\[0.5ex]
   0 & \frac{\partial \hat{f}_1}{\partial \dot{x}_1} & \frac{\partial \hat{f}_1}{\partial \dot{x}_2} & F_{1,3} & I & 0 \\[0.5ex]
	    0 &\frac{\partial \hat{f}_3}{\partial \dot{x}_1} & \frac{\partial \hat{f}_3}{\partial \dot{x}_2} & 0 &F_{3,1}&F_{3,2}  
\end{pmatrix}
\end{eqnarray}
is 1-full with respect to the first column.

We distinguish two cases:
\begin{enumerate}
\item If $F_{3,2} \neq \left( 0, \ldots, 0 \right)$,
  then obviously \eqref{eq:ReducedMatrixInd2OptimalSol}
	is 1-full with respect to the first column if and only if the bordered matrix
	\[
\begin{pmatrix}
	F_{1,3} & I & 0  \\[0.5ex]
	0& \frac{\partial \hat{f}_2}{\partial \dot{x}_1} & \frac{\partial \hat{f}_2}{\partial \dot{x}_2} \\[0.5ex]
	0 & F_{3,1} &F_{3,2} 
\end{pmatrix}=
\left(
\begin{array}{c:cc}
	F_{1,3} & I & 0  \\[0.5ex]
	0& \frac{\partial f_2}{\partial x_1} & \frac{\partial f_2}{\partial x_2} \\[0.5ex] \hdashline
	0 & F_{3,1} &F_{3,2} 
\end{array}
\right)
\]
is non-singular.
Thus the index is two, if and only if
\begin{eqnarray*}
\begin{pmatrix}
	 \frac{\partial f_2}{\partial x_1} \cdot F_{1,3}& \frac{\partial f_2}{\partial x_2} \\[0.5ex]
	 F_{3,1}\cdot F_{1,3} & F_{3,2} 
\end{pmatrix}
\end{eqnarray*}
is non-singular, i.e., iff
\begin{equation} \label{optimal-index2}
F_{3,1}\cdot F_{1,3} - F_{3,2} \cdot
\left( \frac{\partial f_2}{\partial x_2} \right)^{-1} \cdot  \frac{\partial f_2}{\partial x_1} \cdot F_{1,3} \neq 0.
\end{equation}
\item
  If $F_{3,2} = \left( 0, \ldots, 0 \right)$, then the index was one for
  $F_{3,1}\cdot F_{1,3} \neq 0$. 
  Hence, we only have to consider
  $F_{3,2} = \left( 0, \ldots, 0 \right)$ and $F_{3,1}\cdot F_{1,3} = 0$.
  In this case, the matrix \eqref{eq:ReducedMatrixInd2OptimalSol} reads as
\begin{eqnarray*}
\begin{pmatrix}
   F_{1,3} & I & 0 & 0\\[0.5ex]
   0 &F_{3,1}&0 & 0 \\[0.5ex]
	 0 &\frac{\partial f_2}{\partial x_1} & \frac{\partial f_2}{\partial x_2}& 0& 0 &0 & \\[0.5ex]
   0 & \frac{\partial f_1}{\partial x_1} & \frac{\partial f_1}{\partial x_2} & F_{1,3} & I & 0 \\[0.5ex]
	    0 & \frac{\partial \hat{f}_3}{\partial \dot{x}_1} & 0 & 0 &F_{3,1}  & 0
\end{pmatrix}
\end{eqnarray*}
Due to $F_{3,1}\cdot F_{1,3} = 0$, the index is 2 if the rectangular matrix
\begin{eqnarray*}
\begin{pmatrix}
	\frac{\partial f_2}{\partial x_1}\cdot F_{1,3}  & \frac{\partial f_2}{\partial x_2}& 0& 0 & \\[0.5ex]
	\frac{\partial f_1}{\partial x_1} \cdot F_{1,3}  & \frac{\partial f_1}{\partial x_2} & F_{1,3} & I  \\[0.5ex]
	\frac{\partial \hat{f}_3}{\partial \dot{x}_1} \cdot F_{1,3} & 0 & 0 &F_{3,1} 
\end{pmatrix}
\end{eqnarray*}
is 1-full with respect to the first column. Otherwise, the index may be higher or not defined.
\end{enumerate}
In the index-2 case, we also showed that the projector~$T$ is
equal to~(\ref{projector-T}) again.
From a structural point of view, it is of special interest that the projector $Q$ depends on the solution, but in fact only on $Ux$ for the projector $U:=I-T$, i.e., not on the higher-index component.
Thus we obtain $Q(t,Ux)$.

\subsection{Numerical index-computation}
\label{sec:index-comp}

The index of DAEs can also be computed numerically. To this end, we compute all derivatives
with automatic differentiation and check the 1-fullness on
the sequence of matrices
${\cal B}^{[1]}$, ${\cal B}^{[2]}, \ldots $, see~\cite{schwarz-lamour}. 

Since for nonlinear DAEs the index is defined locally, we will obtain
an index statement for a particular consistent initial value.
The algorithm from~\cite{schwarz-lamour} starts with an user-given guess,
computes a corresponding consistent initial value $x_0$ 
and determines the index for that $x_0$. 

Computing the index, several rank decisions have to be made, 
which may be difficult to realise. Recall that above we
 pointed out that in some cases the
 scalar condition \eqref{ftimesf} may determine the index.
 However, if the complete matrices
 ${\cal B}^{[k]}$ contain singular values of different orders of magnitude,
 then the numerically computed index may not be robust.



\section{Index-Analysis}
\label{sec:index}
We identify the index of DAEs from the method of lines
under certain assumptions,
which have an interpretation with respect to the multivariate solution.

\subsection{Analysis for phase condition}
Using the phase condition~(\ref{phasecondition}) for
the MPDAEs~(\ref{mpdae}),
we require the following property in the semi-explicit case.
\begin{condition} \label{cond:phase}
  A solution of the initial-boundary value problem
  (\ref{mpdae-semi-expl}), (\ref{ibvp-semi-expl})
  exists, which satisfies the phase condition~(\ref{phasecondition})
  for the $\ell$th component of the differential variables
  ${ \hat{y}} = (\hat{y}_1,\ldots,\hat{y}_{n_y})^\top$
  together with the property
  $$ \frac{\partial \hat{y}_{\ell}}{\partial t_2}(t_1,0) \neq 0
  \qquad \mbox{for all}\;\; t_1 . $$
\end{condition}
Condition~\ref{cond:phase} is also necessary for the applicability of
Newton's method in a full discretisation of the problem.
We achieve the following conclusion for the method of lines using
a phase condition~(\ref{phase-semi-expl}).

\clearpage

\begin{theorem} \label{thm:phase}
  Let the DAE~(\ref{semi-expl-dae}) be of index one.
  Let the method of lines be convergent.
  If the phase condition~(\ref{phasecondition}) is selected in the
  $\ell$th differential variable satisfying Condition~\ref{cond:phase}
  and the step size~$h$ is sufficiently small,  
  then the differential index of the DAE system
  (\ref{mol-semi-expl}), (\ref{phase-semi-expl})
  in the method of lines is equal to two.
\end{theorem}

\underline{Proof:}

  A differentiation of the differential part in~(\ref{mol-semi-expl}) yields
  \begin{equation} \label{diffdiff}
  {\ddot{ {\bar{y}} }}_i = \textstyle
  \frac{\partial f}{\partial t} +
  \frac{\partial f}{\partial y} {\dot{ {\bar{y}} }}_i +
  \frac{\partial f}{\partial z} {\dot{ {\bar{z}} }}_i -
  \dot{\nu} {\rm D}_i ({\bar{y}}) -
  \nu \textstyle \frac{{\rm d}}{{\rm d}t} {\rm D}_i ({\bar{y}}) .
  \end{equation}
  The term ${\dot{ { \bar{ y}} }}_i$ as well as the
  part $\frac{{\rm d}}{{\rm d}t} {\rm D}_i ({ \bar{y}})$
  can be replaced by the right-hand side of~(\ref{mol-semi-expl}).
  If the DAE~(\ref{semi-expl-dae}) features index one, then the derivatives
  $\dot{ { \bar{ z}}}$ can be eliminated after one differentiation
  of~(\ref{mol-semi-expl}).

  Two differentiations of the phase condition~(\ref{phase-semi-expl}) yield
  $\ddot{\bar{y}}_{1,\ell} = \ddot{\eta}$
  assuming a sufficiently smooth function.
  The function $\eta$ is predetermined and thus not a part of the solution. 
  Now Eq.~(\ref{diffdiff}) implies
  $$ \dot{\nu} \, {\rm D}_{1,\ell} ({ \bar{y}}) =
  r(t,{ \bar{y}},{ \bar{z}},\nu,
  {\ddot{\eta}}) $$
  with some function~$r$.
  It holds that
  $$ \lim_{h \rightarrow 0} \; ({\rm D}_{1,\ell} ({ \bar{y}})) (t_1) =
  \textstyle \frac{\partial \hat{y}_{\ell}}{\partial t_2}(t_1,0)
  \qquad \mbox{for each} \;\; t_1 , $$
  because the method of lines is assumed to be convergent.
  Since a compact interval $t_1 \in [0,t_{\rm end}]$ is assumed,
  even uniform convergence is given. 
  If the step size~$h$ is sufficiently small, then the finite difference
  approximation is also non-zero due to Condition~\ref{cond:phase}.
  We obtain an ODE for $\nu$ by a division.
  It follows that the index is exactly two.
\hfill $\Box$

\medskip

If we select the phase condition~(\ref{phasecondition}) in an
algebraic variable, then the above derivation is not feasible any more.
It follows that the differential index is at least two.

The above derivations are in agreement to the structural analysis in
Section~\ref{sec:structure-phase}.
We concluded that the index must be larger than one
in Section~\ref{sec:structure-phase-index1}.
The proof of Theorem~\ref{thm:phase} shows the property~(\ref{ftimesf}),
which implies that the index is equal to two.

Considering models by (implicit) ODEs, the phase condition is naturally
chosen in a differential variable.
We obtain directly the following statement.
\begin{corollary} \label{cor:phase}
  If the system~(\ref{dae}) consists of ODEs (index zero),
  then the system~(\ref{mol}), (\ref{phase-mol}) from the method of lines
  exhibits the differential index two provided that
  the assumptions of Theorem~\ref{thm:phase} are satisfied. 
\end{corollary}

\subsection{Analysis for optimal solutions}
We consider the functional~(\ref{functional}) of the optimisation
and the associated necessary condition~(\ref{opt-disc2}) now.
The following property is required.
\begin{condition} \label{cond:non-constant}
  The initial-boundary value problem
  (\ref{mpdae-semi-expl}), (\ref{ibvp-semi-expl}) has a solution,
  where at least one variable involved in the optimisation is
  non-constant in the fast time scale, i.e.,
  \begin{equation} \label{non-constant}
    \int_0^1 \left\| { W}_y^{\frac{1}{2}}
    {\textstyle \frac{\partial { \hat{y}}}{\partial t_2}}
    \right\|_2^2 \; {\rm d}t_2 +
    \int_0^1 \left\| { W}_z^{\frac{1}{2}}
    {\textstyle \frac{\partial { \hat{z}}}{\partial t_2}}
    \right\|_2^2 \; {\rm d}t_2 > 0
    \qquad \mbox{for all} \;\; t_1 
  \end{equation}
with
${ W}_y = {\rm diag} ( w_1^y , \ldots , w_{n_y}^y )$ as well as
${ W}_z = {\rm diag} ( w_1^z , \ldots , w_{n_z}^z )$
and the Euclidean norm~$\| \cdot \|_2$.
\end{condition}
If there is such a non-constant solution, then all solutions of the
initial-boundary value problem satisfy this property due to
a transformation, cf.~\cite{pulch08proc}.
Condition~\ref{cond:non-constant} is given in most of the cases
provided that the weights do not all vanish.
$W_y = 0$ or $W_z = 0$ imply that the second term or the first term
is positive, respectively.
Furthermore, conditions of this type were also required for the existence and
uniqueness of optimal solutions in~\cite{kugelmann-pulch}.

\begin{theorem} \label{thm:opt}
  Let the DAE~(\ref{semi-expl-dae}) have index one.
  Let the method of lines be convergent and a sufficiently small
  step size~$h$ be given.
  The differential index of the system
  (\ref{mol-semi-expl}), (\ref{opt-disc2}) from the method of lines
  exhibits the following relations provided that
  Condition~\ref{cond:non-constant} is satisfied.
  \begin{itemize}
  \item[i)] If $w_1^z = \cdots = w_{n_z}^z=0$, then the index is one.
  \item[ii)] If $w_1^z + \cdots + w_{n_z}^z > 0$, then the index is at
    least two. \\
    The index is exactly two for linear algebraic constraints.
  \end{itemize}
\end{theorem}

\underline{Proof:}

Case~(i):
A differentiation of the necessary constraint~(\ref{opt-disc2}) yields
\begin{equation} \label{opt-disc3}
\begin{array}{l}
\displaystyle \sum_{i=1}^m \bigg( \sum_{\ell = 1}^{n_y} 
w_\ell^y \left( \ddot{\bar{y}}_{i,\ell} \cdot
{\rm D}_{i,\ell} ({ \bar{y}}) +
\dot{\bar{y}}_{i,\ell} \cdot
{\textstyle \frac{\rm d}{{\rm d}t}} {\rm D}_{i,\ell} ({ \bar{y}}) \right) \\
\mbox{} \qquad + \displaystyle
\sum_{\ell = 1}^{n_z} 
w_\ell^z \left( \ddot{\bar{z}}_{i,\ell} \cdot
{\rm D}_{i,\ell} ({ \bar{z}}) +
\dot{\bar{z}}_{i,\ell} \cdot
{\textstyle \frac{\rm d}{{\rm d}t}} {\rm D}_{i,\ell} ({ \bar{z}}) \right)
\bigg) = 0 . \\
\end{array} 
\end{equation}
The second inner sum vanishes due to $w_{\ell}^z = 0$ for all~$\ell$.
One differentiation of the differential part in~(\ref{mol-semi-expl})
allows to insert~(\ref{diffdiff}) into~(\ref{opt-disc3}).
It follows that
$$ \sum_{i=1}^m \sum_{\ell = 1}^{n_y} 
  w_\ell^y \left( (- \dot{\nu} {\rm D}_{i,\ell} ({ \bar{y}}) +
  r_{\ell}(t,{ \bar{y}},\dot{\bar{y}},\dot{\bar{z}},\nu) ) \cdot
  {\rm D}_{i,\ell} ({ \bar{y}}) + \dot{\bar{y}}_{i,\ell} \cdot
  {\textstyle \frac{\rm d}{{\rm d}t}} {\rm D}_{i,\ell} ({ \bar{y}})
  \right) = 0 $$
with functions $r_1,\ldots,r_{n_y}$.
We obtain 
\begin{equation} \label{nudot}
  \dot{\nu} \; h \sum_{i=1}^m \sum_{\ell = 1}^{n_y}
  w_{\ell}^y \left( {\rm D}_{i,\ell} ({ \bar{y}}) \right)^2 =
  h \, \tilde{r} (t,{ \bar{y}},\dot{\bar{y}},\dot{\bar{z}},\nu) 
\end{equation}
with another function~$\tilde{r}$.
It holds that
\begin{equation} \label{limit-integral}
  \lim_{h \rightarrow 0}
  h \sum_{i=1}^m \sum_{\ell = 1}^{n_y}
  w_{\ell}^y \left( {\rm D}_{i,\ell} ({ \bar{y}}) \right)^2 =
  \int_0^1 \left\| { W}_y^{\frac{1}{2}}
      {\textstyle \frac{\partial { \hat{y}}}{\partial t_2}}
      \right\|_2^2 \; {\rm d}t_2
\end{equation}
due to the assumption of a convergent method of lines.
If Condition~\ref{cond:non-constant} and thus~(\ref{non-constant})
is satisfied, then the right-hand side of~(\ref{limit-integral}) is
positive.
Hence a non-zero term appears for sufficiently small step size.
In Eq.~(\ref{nudot}), a division yields an equality for~$\dot{\nu}$.
Since we required just a single differentiation,
the differential index is one.

Case~(ii):
Since algebraic variables are involved in the optimisation,
a derivation as in case~(i) is not feasible.
Thus the differential index becomes larger than one.

Let the algebraic constraints of the system~(\ref{semi-expl-dae})
be linear, i.e., 
$$ G_y y + G_z z + b_g = 0 $$
with constant matrices $G_y,G_z$ and a vector-valued function $b_g$.
The index-1 property guarantees a non-singular matrix~$G_z$. 
We show that the condition~(\ref{optimal-index2}) from the structural
analysis is satisfied, which implies a differential index equal to two.
It holds that $x_1 = \bar{y}$ and $x_2 = \bar{z}$.

The difference formula~(\ref{difference-formula}) is assumed to be
convergent.
It follows that the sum of all coefficients~$\alpha_j$ is equal to zero.
We calculate
\begin{equation} \label{diff-change}
  \begin{array}{rcl}
    G_z^{-1} G_y \; {\rm D}_i ({ \bar{y}}) & = & 
   \displaystyle \frac{1}{h} \displaystyle \sum_{j=-q}^p \alpha_j
   G_z^{-1} G_y \;
   \bar{y}_{i+j} \;\; = \;\;  
        \displaystyle \frac{1}{h} \displaystyle \sum_{j=-q}^p \alpha_j
        \left( - \bar{z}_{i+j} - G_z^{-1} b_g \right) \\
        & = &
        \displaystyle - \frac{1}{h} \displaystyle
        \sum_{j=-q}^p \alpha_j \bar{z}_{i+j} -
        \displaystyle \frac{1}{h} G_z^{-1} b_g
        \displaystyle \sum_{j=-q}^p \alpha_j
         \;\; = \;\;  
        - {\rm D}_i ({ \bar{z}}) \\
  \end{array}
\end{equation}
for $i=1,\ldots,m$.
The application of the difference operator can be described by
${\rm D} ({ \bar{y}}) = S_1 \bar{y}$ and
${\rm D} ({ \bar{z}}) = S_2 \bar{z}$
with constant matrices $S_1 \in \real^{mn_y \times mn_y}$ and
$S_2 \in \real^{mn_z \times mn_z}$.
It follows that $F_{1,3} = - S_1 x_1$.
In~(\ref{optimal-index2}), the Jacobian matrices become
$\frac{\partial f_2}{\partial x_1} = I_m \otimes G_y$ and
$\frac{\partial f_2}{\partial x_2} = I_m \otimes G_z$
using the Kronecker product and the identity matrix
$I_m \in \real^{m \times m}$.
Hence these matrices are constant and block-diagonal.
Due to~(\ref{diff-change}), we obtain
\begin{equation} \label{diff-change2}
  \left( \textstyle \frac{\partial f_2}{\partial x_2} \right)^{-1} \;
  \textstyle \frac{\partial f_2}{\partial x_1} \; S_1 x_1 = 
  ( I_m \otimes (G_z^{-1} G_y) ) S_1 x_1
  = - S_2 x_2 .
\end{equation}
Let $W_1 := I_m \otimes W_y$ and $W_2 := I_m \otimes W_z$.
Hence $W_1,W_2$ are positive semi-definite.
Concerning the necessary condition~(\ref{opt-disc2})
written in the form~(\ref{eq:opt-general}),
it follows that
$$ F_{3,k} = (W_k S_k x_k)^\top \qquad \mbox{for} \;\; k=1,2 . $$
Now we investigate the property~(\ref{optimal-index2}).
Eq.~(\ref{diff-change2}) yields
$$ \begin{array}{cl}
 & F_{3,1} F_{1,3} - F_{3,2} 
\left( \frac{\partial f_2}{\partial x_2} \right)^{-1}
\frac{\partial f_2}{\partial x_1} F_{1,3} \\
= & (W_1 S_1 x_1)^\top (-S_1 x_1) - (W_2 S_2 x_2)^\top
\left( \frac{\partial f_2}{\partial x_2} \right)^{-1} 
\frac{\partial f_2}{\partial x_1} (-S_1 x_1) \\[2ex]
= & - \left[ x_1^\top S_1^\top W_1 S_1 x_1
  + x_2^\top S_2^\top W_2 S_2 x_2 \right] \\[2ex]
= & - \big[ \underbrace{ (S_1 x_1)^\top W_1 (S_1 x_1) }_{ \ge 0}
  + \underbrace{ (S_2 x_2)^\top W_2 (S_2 x_2) }_{\ge 0} \big] . \\
\end{array} $$
Condition~\ref{cond:non-constant} guarantees that one of the two
non-negative terms is positive for sufficiently small step size.
Thus the property~(\ref{optimal-index2}) is satisfied.
\hfill $\Box$

\medskip

Again the above results are in agreement to the structural analysis
in Section~\ref{sec:structural-opt}.
In case~(i), the proof of Theorem~\ref{thm:opt} shows that
the property~(\ref{ftimesf}) is fulfilled,
which guarantees an index-one system.

A system of (implicit) ODEs is equivalent to
a system~(\ref{semi-expl-dae}) without an algebraic part.
At least one weight of the minimisation is positive and thus
case~(i) can be applied.

\begin{corollary} \label{cor:opt}
  If the system~(\ref{dae}) consists of ODEs (index zero),
  then the system~(\ref{mol}), (\ref{opt-disc}) from the method of lines
  exhibits the differential index one
  provided that the assumptions of Theorem~\ref{thm:opt} 
  are satisfied.
\end{corollary}


\section{Illustrative Example}
\label{sec:example}
We show results of numerical simulations,
where the systems of DAEs following from the method of lines are solved.

\begin{figure}[t]
\begin{center}
\includegraphics[width=14cm]{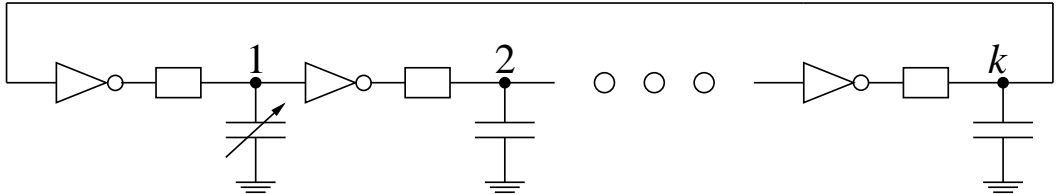}
\end{center}
\caption{Electric circuit of a $k$-stage ring oscillator.}
\label{fig:circuit}
\end{figure}

\subsection{Modelling and simulation of a ring oscillator}
\label{sec:ringoscillator}
In~\cite[Sect.IV.B]{lai}, the ODE model of a three-stage ring oscillator
was considered.
We extend this example to a ring oscillator with $k$-stages
depicted in Figure~\ref{fig:circuit}.
An odd number~$k$ is required such that the circuit exhibits the
desired behaviour.
The circuit consists of capacitances, resistances and inverters.
We model this electric circuit by a nonlinear semi-explicit system of
DAEs~(\ref{semi-expl-dae}) with differential index~one:
\begin{equation} \label{ringosc}
  \begin{array}{rclrcl}
    \dot{u}_1 & = & \imath_1 / (C b) , \qquad & 
    0 & = & R \, \imath_1 - (\tanh(Gu_k) - u_1) , \\
    \dot{u}_2 & = & \imath_2 / C , &
    0 & = & R \, \imath_2 - (\tanh(Gu_1) - u_2) , \\
    & \vdots & & & \vdots & \\
    \dot{u}_k & = & \imath_k / C , & 
    0 & = & R \, \imath_k - (\tanh(Gu_{k-1}) - u_k) . \\
  \end{array}
\end{equation}
The unknowns are the node voltages $u_1,\ldots,u_k$ and the
branch currents $\imath_1,\ldots,\imath_k$.
Thus the dimension of the system becomes~$2k$.
The constant parameters $C,R,G$ have to be predetermined.
The first capacitance is controlled by an independent input signal~$b$.
A constant input like $b \equiv 1$ implies an autonomous
system~(\ref{ringosc}) with a stable periodic solution.
A slowly varying input causes both
amplitude modulation and frequency modulation. 

Thus we apply the multidimensional model,
where the DAE system~(\ref{ringosc}) changes into a semi-explicit
MPDAE system~(\ref{mpdae-semi-expl}).
The method of lines yields the system~(\ref{mol-semi-expl}),
where $m=100$ lines are applied in the following.
We always use the BDF formula~(\ref{bdf}) of order two
in the finite differences~(\ref{difference-formula}).
The implicit Euler method (BDF-1) produces the numerical solution of
initial value problems, where always 200 time steps are performed in
a global interval of the slow time scale.

The method of lines requires an additional constraint.
On the one hand, we employ a phase condition from
Section~\ref{sec:phasecondition}.
Choosing the first node voltage, the phase condition~(\ref{phasecondition})
reads as
\begin{equation} \label{phase-ringosc}
  \hat{u}_1(t_1,0) = 0 \qquad \mbox{for all} \;\; t_1 .
\end{equation}
Theorem~\ref{thm:phase} guarantees that the index of the DAE system
is equal to two for sufficiently small step size in the method of lines.
On the other hand, we apply the necessary condition from
Section~\ref{sec:optimalsolution} in the discretised form~(\ref{opt-disc}).
Three choices of weights are investigated:
\begin{enumerate}
\item[a)] $W_y = I_k$, $W_z = I_k$,
\item[b)] $W_y = I_k$, $W_z = 0$,
\item[c)] $W_y = 0$, $W_z = I_k$,
\end{enumerate}
with the identity matrix $I_k \in \real^{k \times k}$.
Theorem~\ref{thm:opt} implies that the index of the DAE system
is equal to one in the case~(b) with sufficiently small step size
and at least two in the cases~(a) and~(c).

We examine the ring oscillator for the two different choices
$k=3$ and $k=11$ of the inverter number.
Another numerical simulation of the three-stage ring oscillator
is also reported in~\cite{pulch15}.

\subsection{Simulation of three-stage ring oscillator}
\label{sec:three-stage}
In the system~(\ref{ringosc}), the physical parameters are fixed
to $C = 10^{-6}$, $R = 10^3$, $G=-5$.
We choose the harmonic oscillation
$$ b(t) = 1 + \textstyle \frac{1}{2} \sin \left( \frac{2\pi}{T} t \right) $$
with the period~$T=1$ as input signal.
The global time interval $t \in [0,T]$ is considered in the
numerical simulation.
The dimension of DAE systems becomes $mk+1=601$ in the method of lines.

We used the algorithm from~\cite{schwarz-lamour} outlined
in Section~\ref{sec:index-comp},
which computes the differential index as well as consistent initial values 
numerically.
Therein, an index of one is verified for the optimisation case~(b),
whereas the other optimisation cases and the phase condition
result in an index of two.
However, the algorithm is not able to determine the index for
more critical parameters $C,R,G$,
because the matrices become ill-conditioned in the rank decisions.

\begin{figure}
\begin{center}
\includegraphics[width=6.5cm]{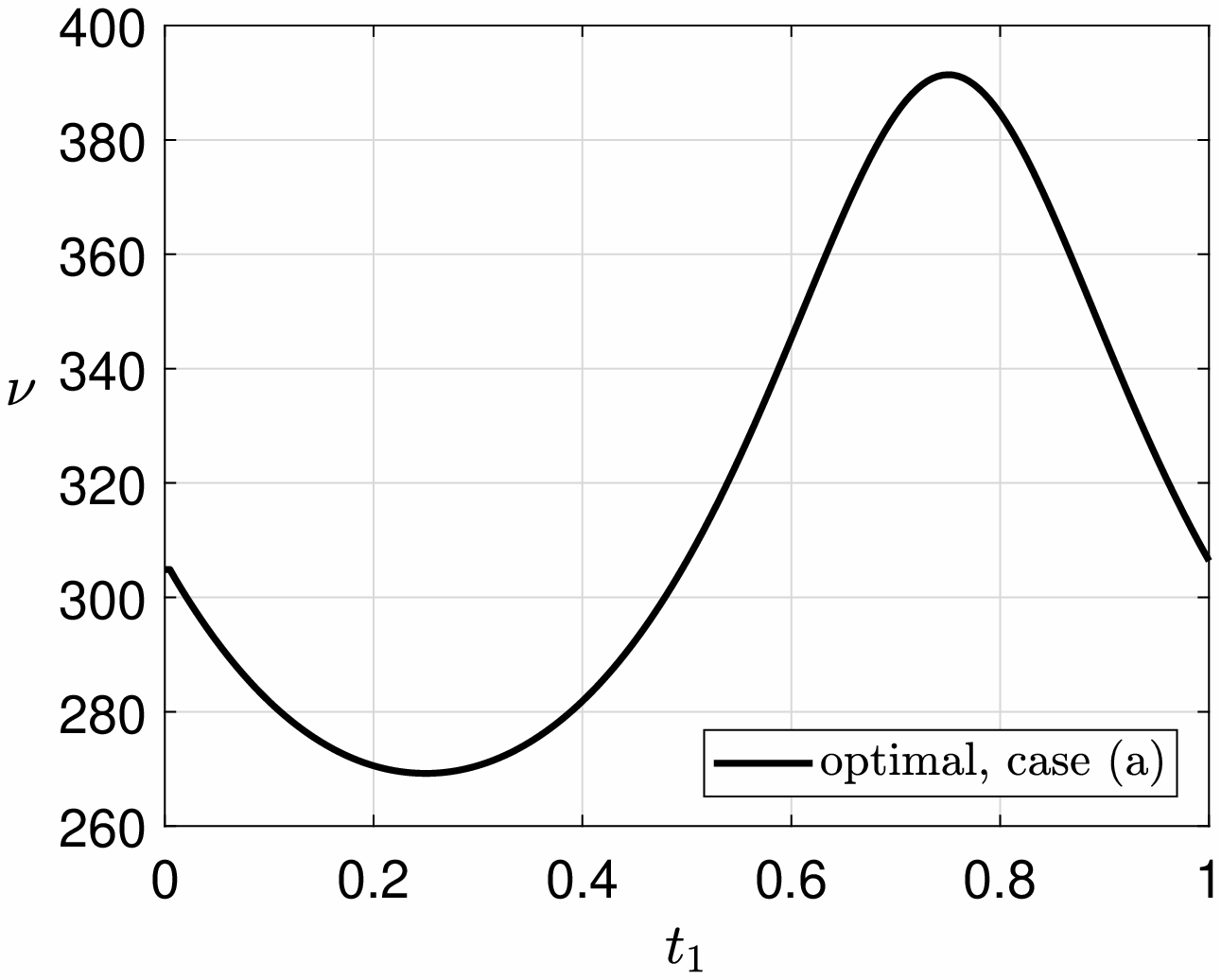}
\hspace{5mm}
\includegraphics[width=6.5cm]{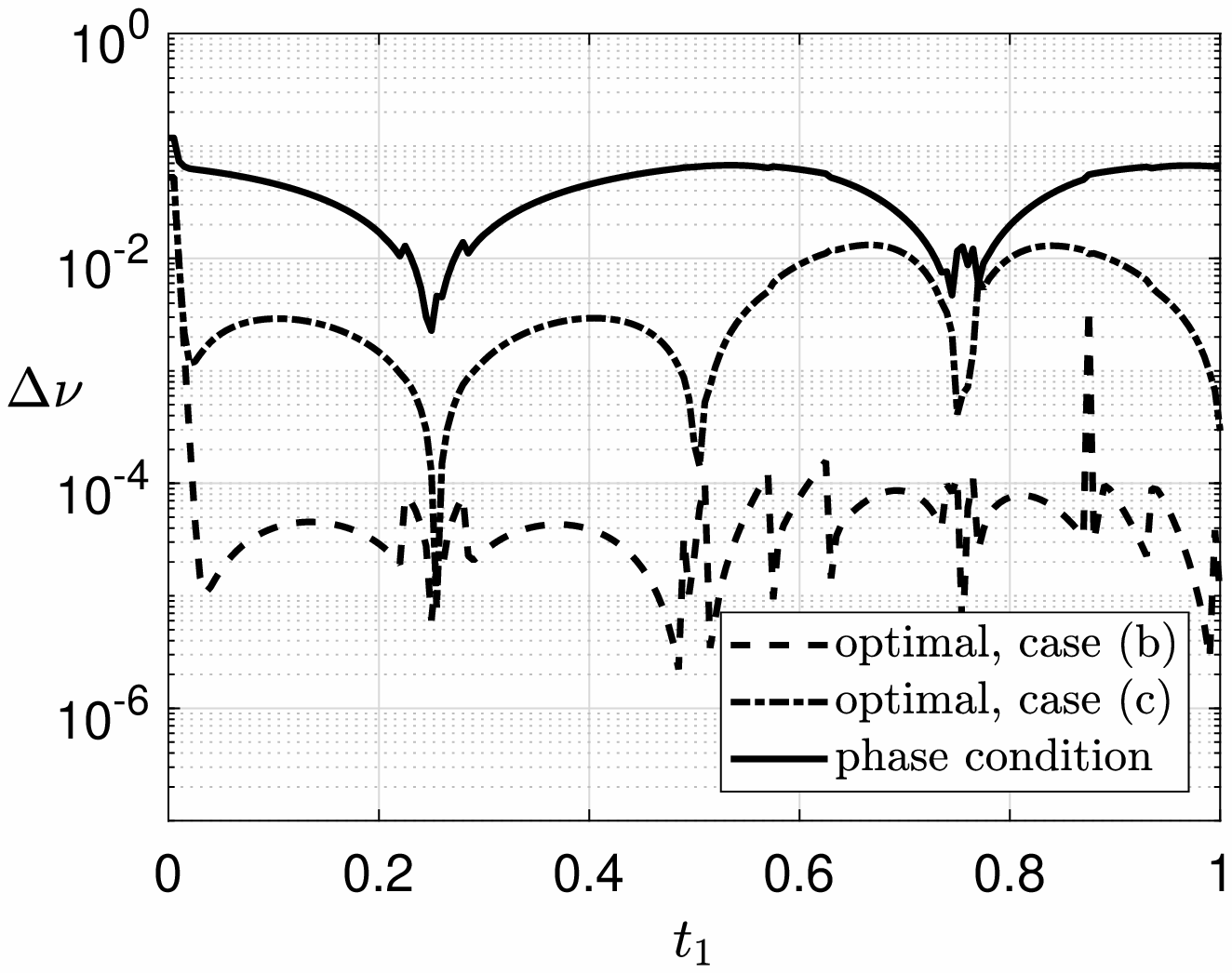}
\end{center}
\caption{Local frequency function for the optimal solution in case~(a) (left)
  and differences to the local frequency functions of the other cases (right)
  in three-stage ring oscillator.}
\label{fig:three-frequency}
\end{figure}

\begin{figure}
  \begin{center}
    node voltage
    
    \includegraphics[width=6.5cm]{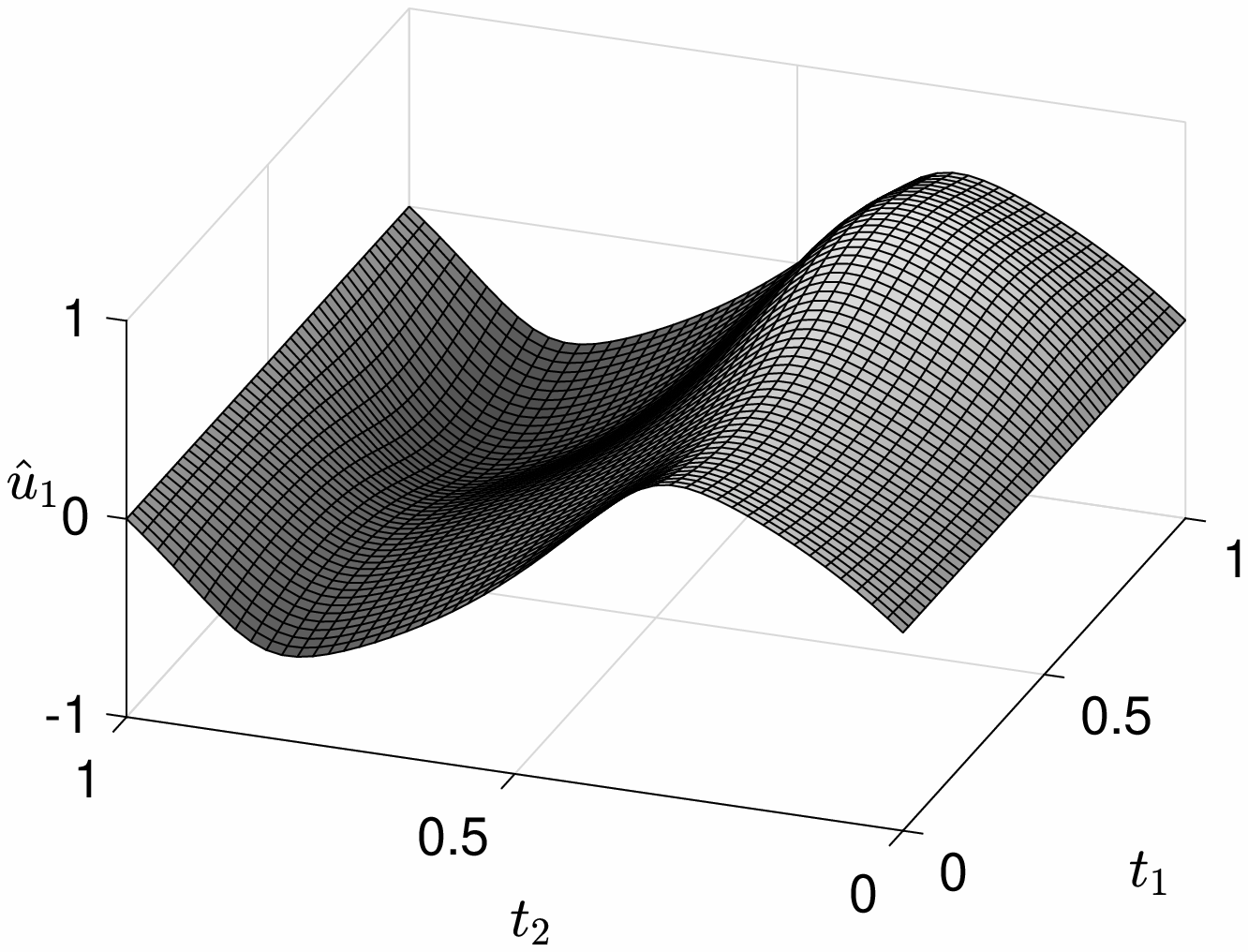}
    \hspace{5mm}
    \includegraphics[width=6.5cm]{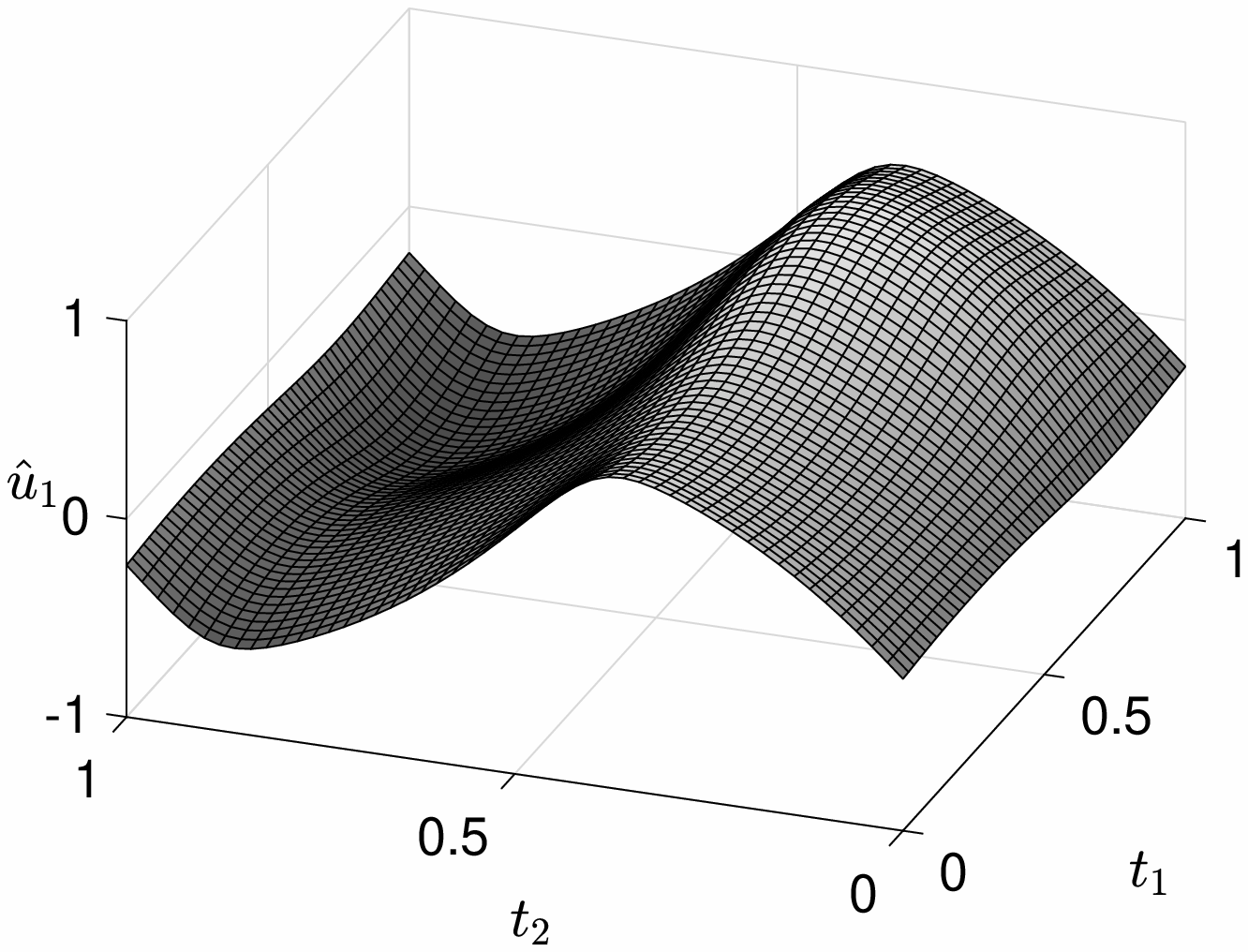}

    branch current
    
    \includegraphics[width=6.5cm]{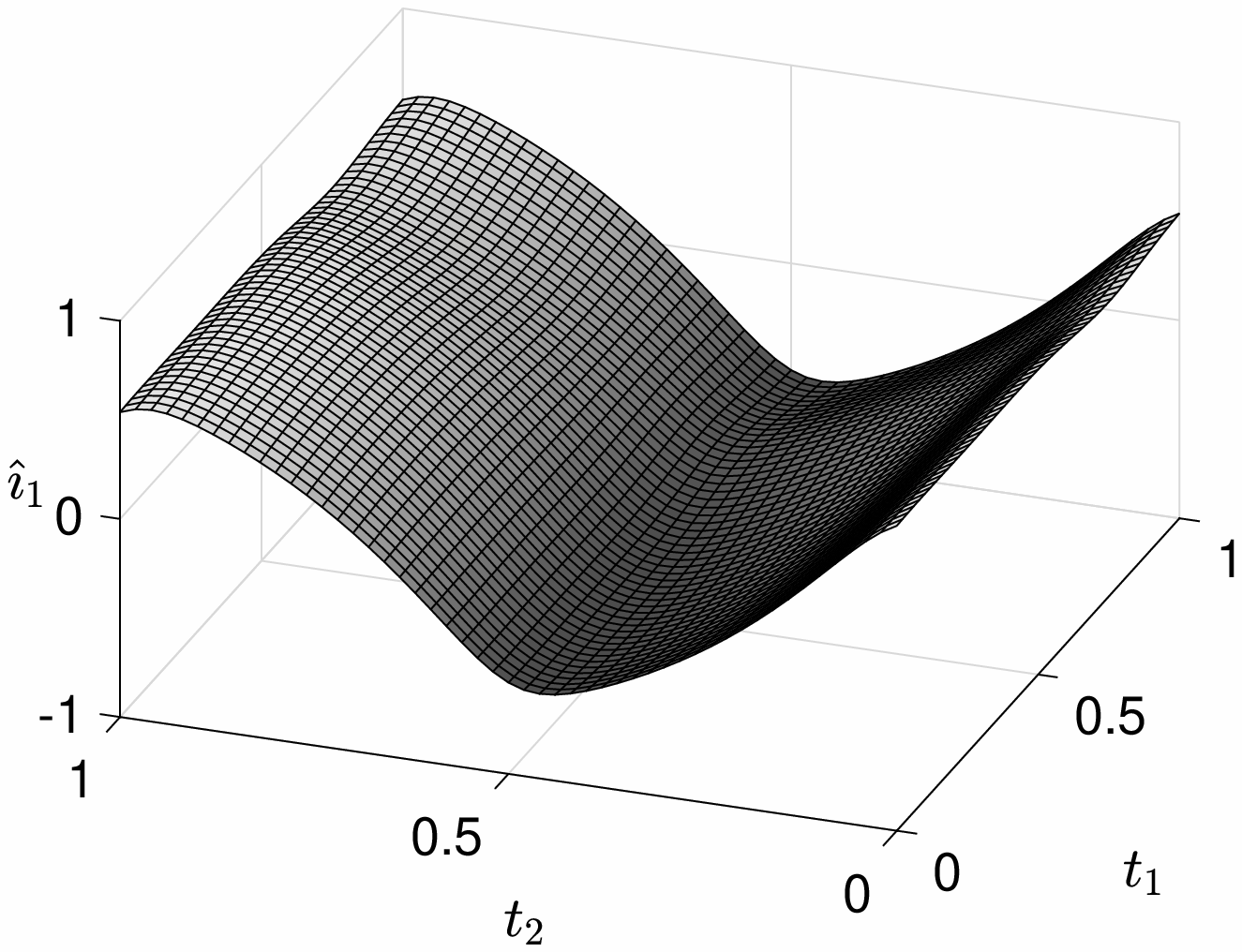}
    \hspace{5mm}
    \includegraphics[width=6.5cm]{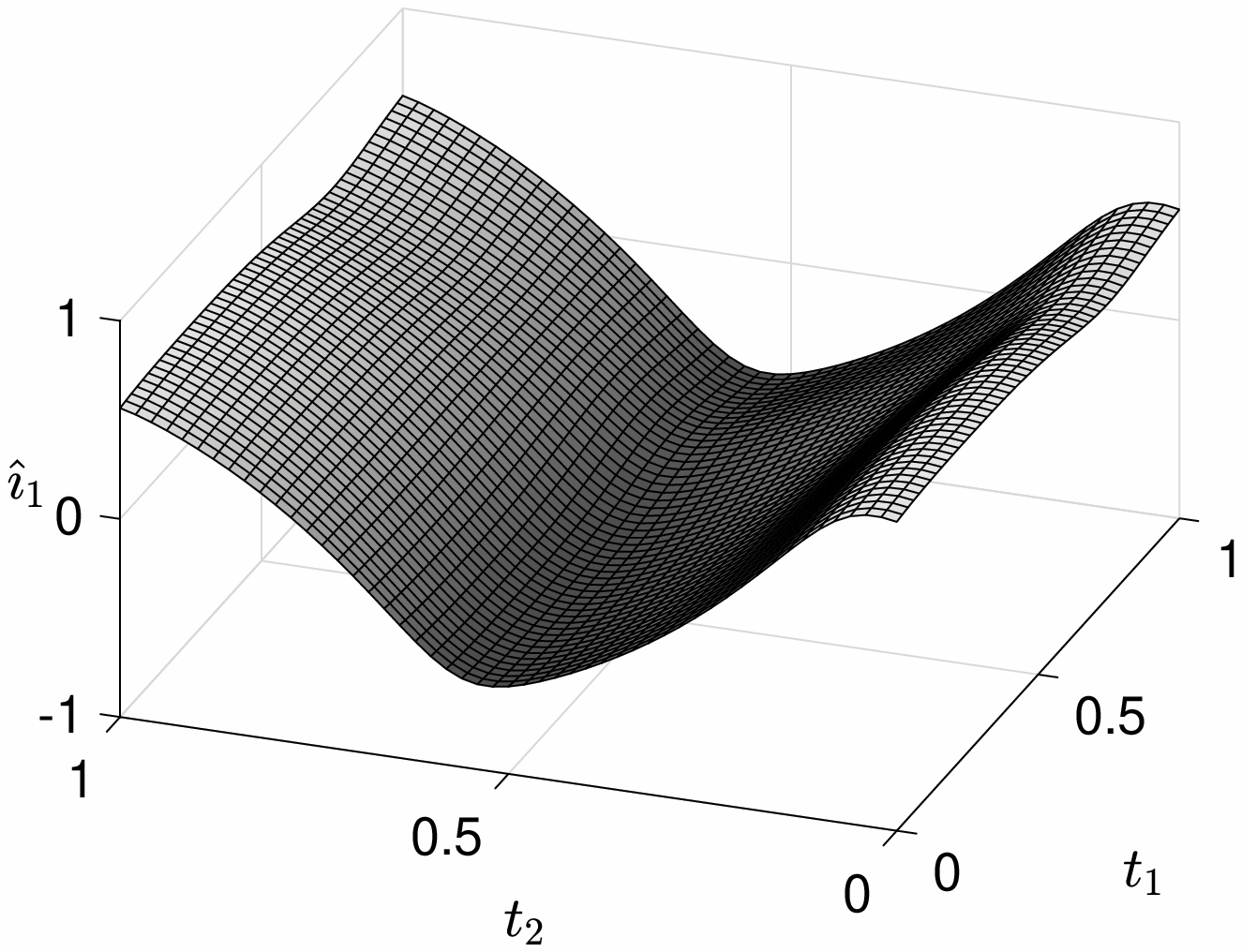}
  \end{center}
\caption{Numerical solutions from method of lines with phase condition (left)
  and optimality condition (right) for node voltage as well as branch current
  in three-stage ring oscillator.}
\label{fig:three-mvfs}
\end{figure}

Now initial value problems of the DAE systems are solved numerically,
as described in Section~\ref{sec:ringoscillator},
using the computed consistent initial values.
Figure~\ref{fig:three-frequency} (left) shows the local frequency
function for the optimal solution in case~(a).
The modulus of the differences between the other local frequencies
and this function are illustrated by Figure~\ref{fig:three-frequency} (right),
where a semi-logarithmic scale is used due to different orders of magnitudes.
The phase condition causes the largest difference. 
Figure~\ref{fig:three-mvfs} depicts the resulting MPDAE solutions
for the first node voltage and the first branch current
with the phase condition and the optimisation case~(a), respectively.
We observe that the functions for the phase condition and
the optimisation are similar.

\subsection{Simulation of eleven-stage ring oscillator}
\label{sec:eleven-stage}

Now the parameters are set to $C = 2 \cdot 10^{-12}$, $R = 10^3$, $G=-5$
in the system~(\ref{ringosc}), which includes a more realistic value
of the capacitance.
A stable periodic solution around 30 MHz emerges for a constant input
$b \equiv 1$.
We supply the input signal 
$$ b(t) = 1 + 2 \textstyle \sin^2 \left( \frac{2\pi}{T} t \right) $$
with the forced time rate $T = 10^{-4}$ (1 kHz),
which causes about 3100 oscillations in the global interval
$t \in [0,T]$.
As initial values for the associated MPDAE system~(\ref{mpdae-semi-expl}),
we take values from an approximation of the stable periodic solution
for $b \equiv 1$.
Now the dimension reads as $mk+1=2201$ in DAE systems from
the method of lines.

\begin{figure}
\begin{center}
\includegraphics[width=6.5cm]{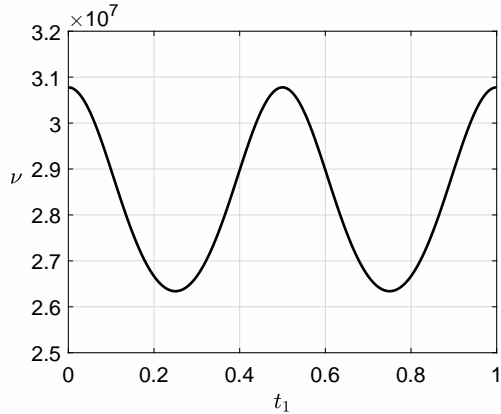}
\end{center}
\caption{Local frequency function for the optimal solution of case~(b) 
  in eleven-stage ring oscillator.}
\label{fig:eleven-frequency}
\end{figure}

\begin{figure}
  \begin{center}
    node voltage
    
    \includegraphics[width=6.5cm]{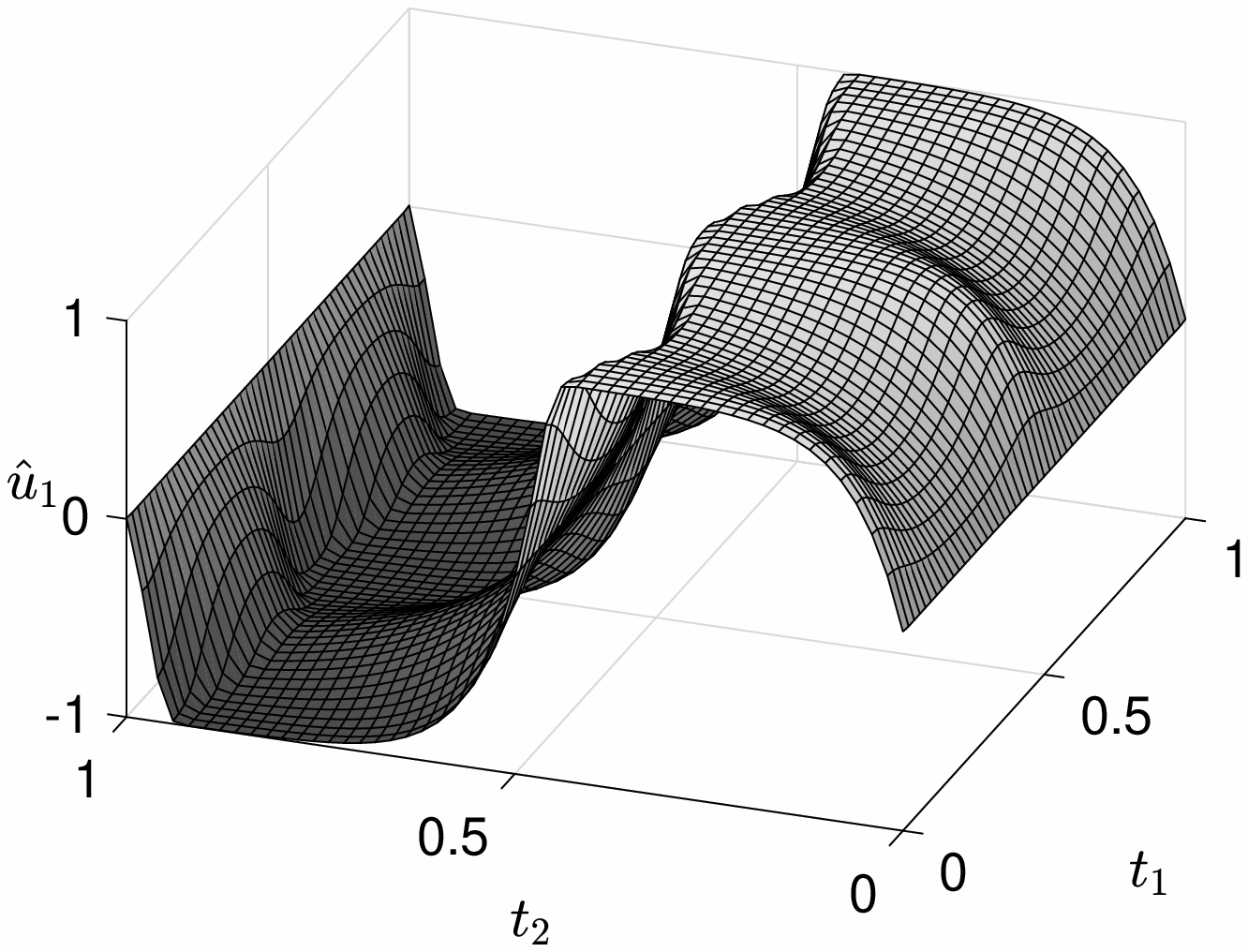}
    \hspace{5mm}
    \includegraphics[width=6.5cm]{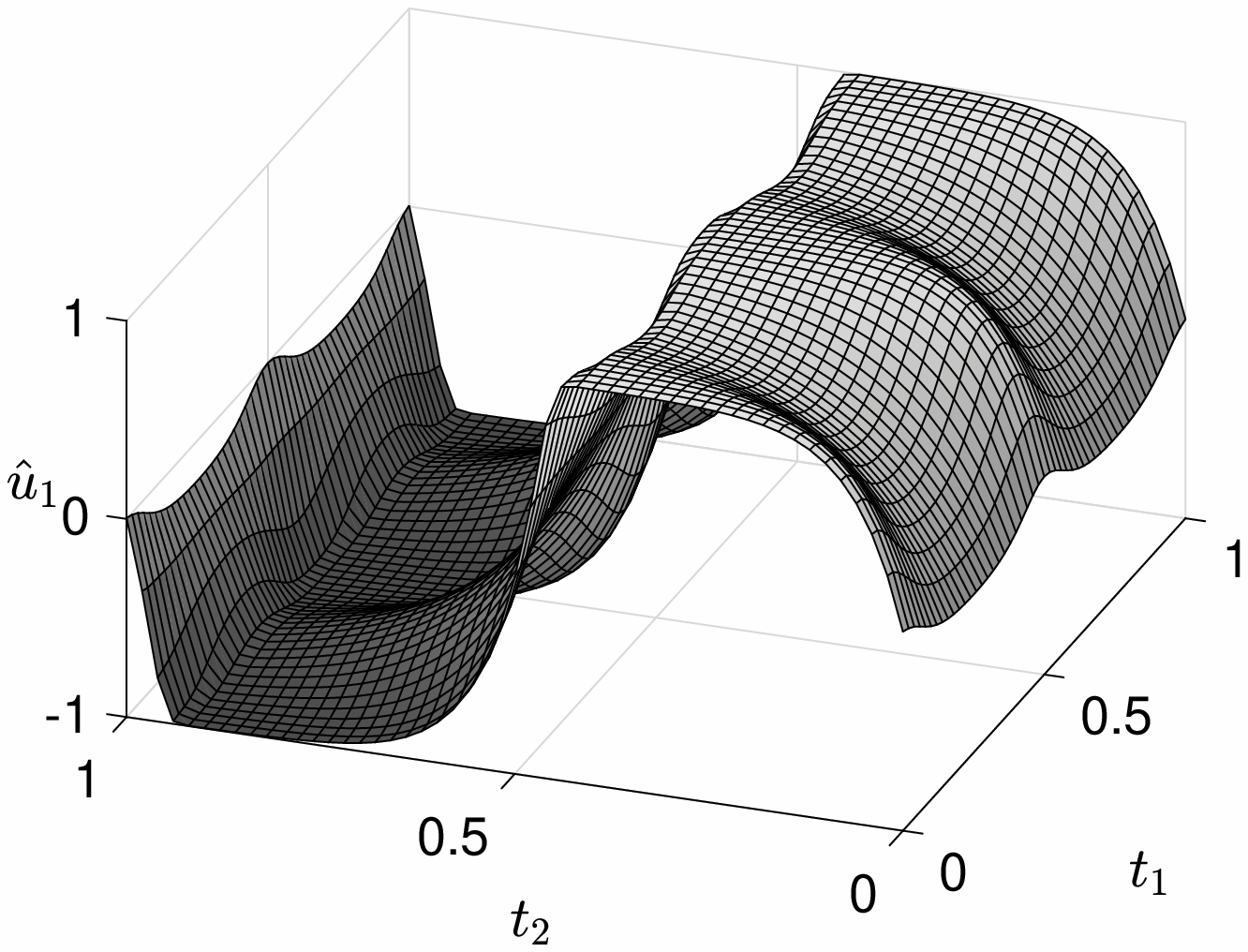}

    branch current
    
    \includegraphics[width=6.5cm]{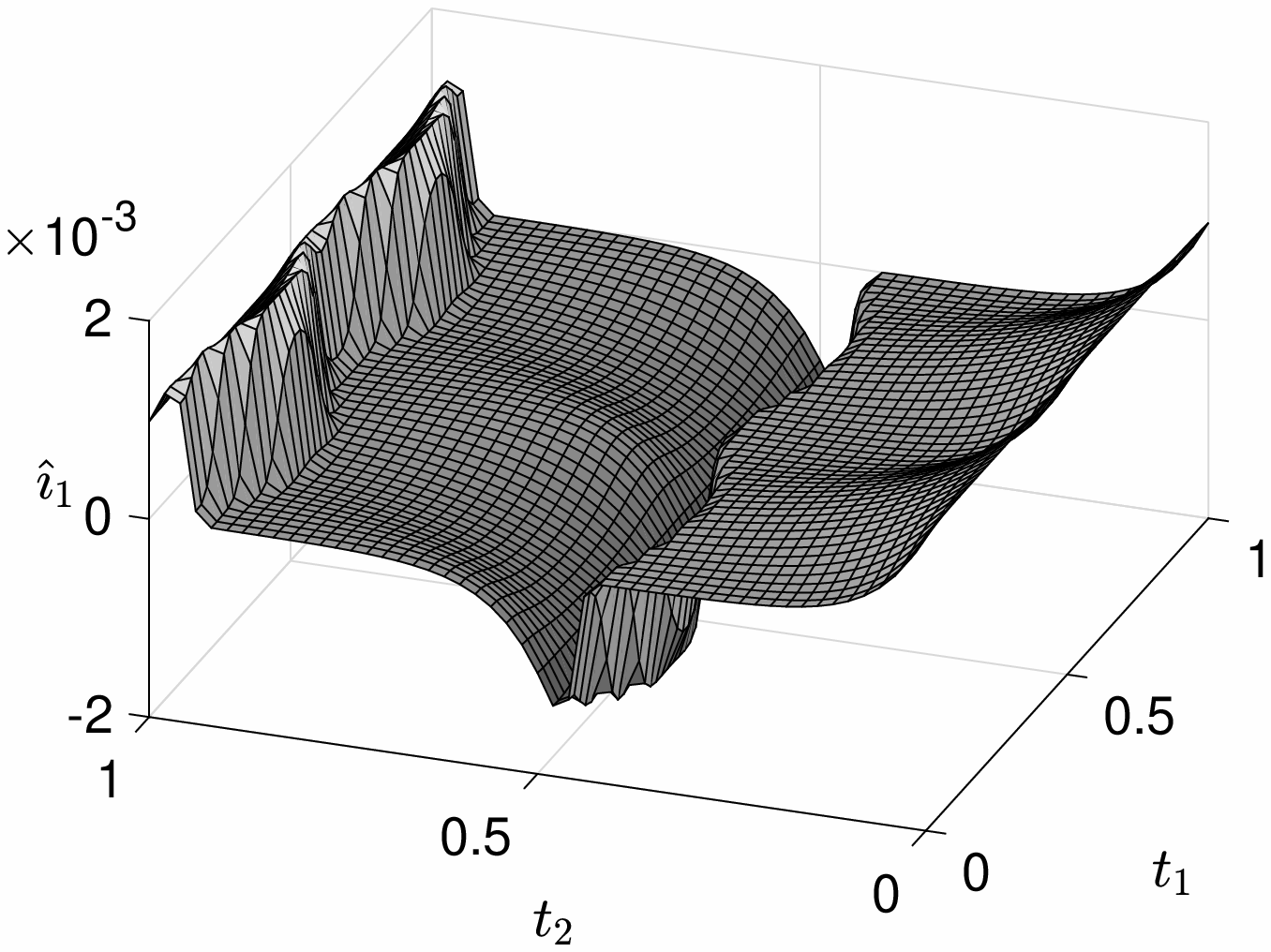}
    \hspace{5mm}
    \includegraphics[width=6.5cm]{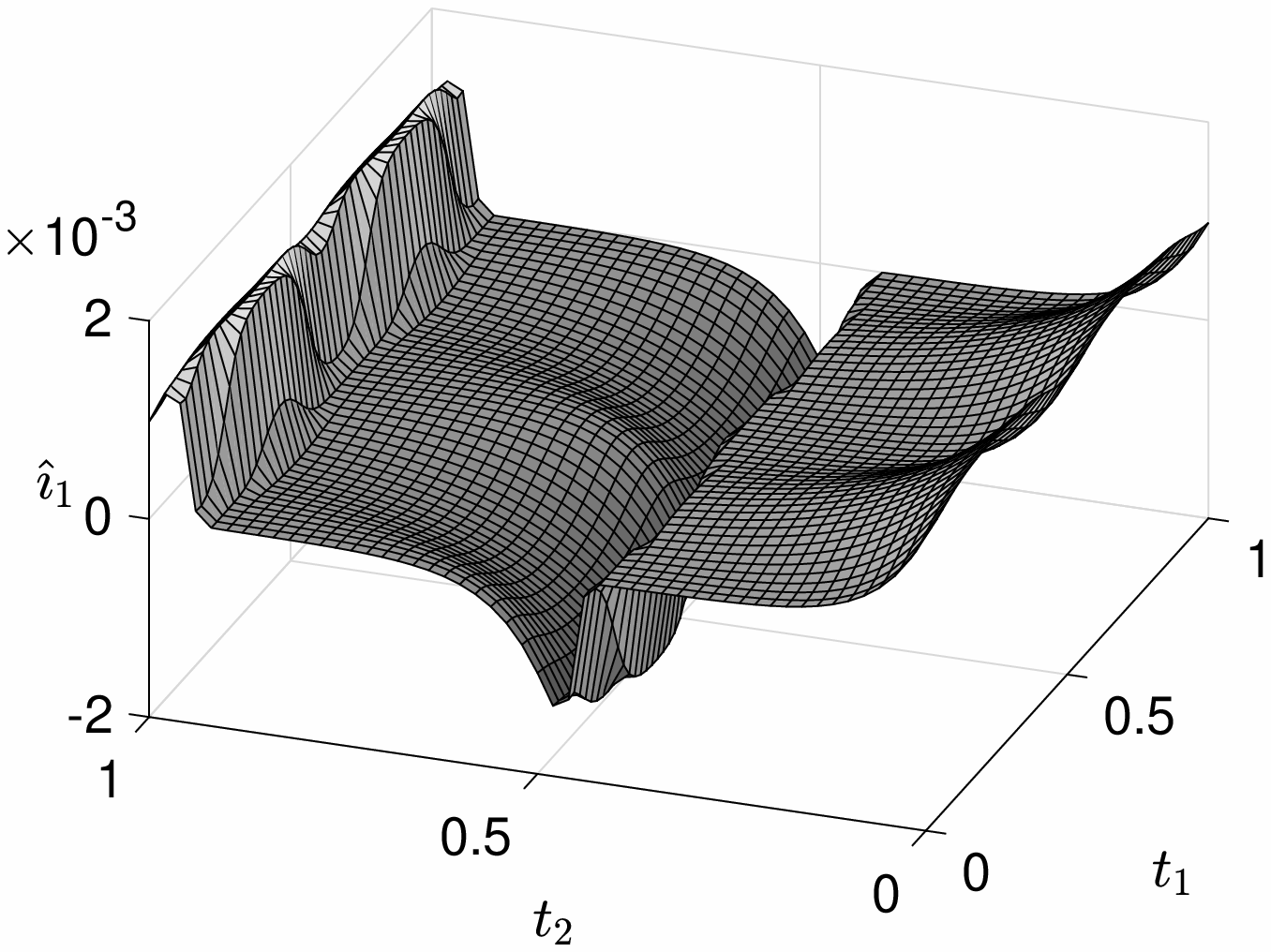}
  \end{center}
\caption{Numerical solutions from method of lines with phase condition (left)
  and optimality condition (right) for node voltage as well as branch current
  in eleven-stage ring oscillator.}
\label{fig:eleven-mvfs}
\end{figure}

We investigate only case~(b) of the optimisation and the case of the
phase condition~(\ref{phase-ringosc}),
where an index-one system and an index-two system,
respectively, is guaranteed for sufficiently small step size
in the method of lines.
We apply the same initial values in both situations.
Thus the initial values are just nearly consistent.
The slow time scale $t_1 \in [0,10^{-4}]$ is standardised to
$t_1 \in [0,1]$ in the plots.
Figure~\ref{fig:eleven-frequency} illustrates the local frequency function
of the optimal solution, whereas the relative difference to the
local frequencies of the phase condition is in a magnitude of just 0.01\%.
The multidimensional solutions of the first node voltage as well as
the first branch current are shown in Figure~\ref{fig:eleven-mvfs}.

\begin{figure}
\begin{center}
  \mbox{} \hspace{1cm}
  \includegraphics[width=12cm]{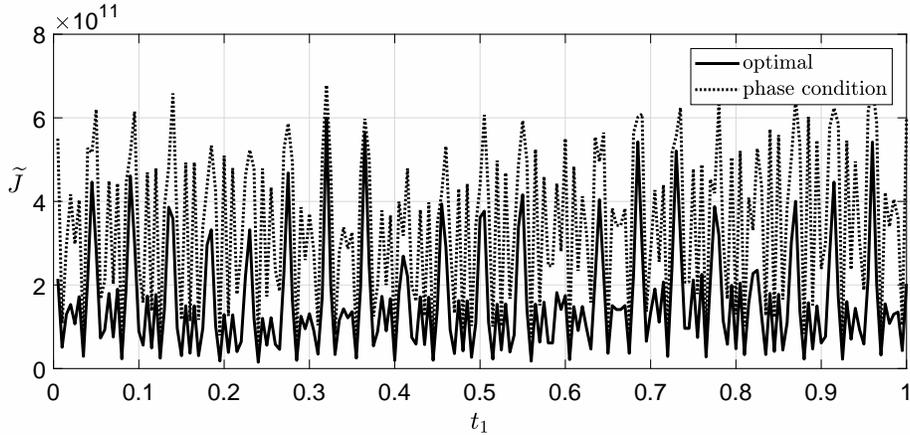}
\end{center}
\caption{Functional~(\ref{functional-pw}) for the optimal solution and
  for the solution with phase condition
  in eleven-stage ring oscillator.}
\label{fig:eleven-functional}
\end{figure}

Finally, we employ the MPDAE solutions to evaluate approximately
the functional~(\ref{functional-pw}) from the minimisation
using the weights of case~(b).
Figure~\ref{fig:eleven-functional} shows these approximations.
Therein, the optimality of the solution obtained by the
necessary condition~(\ref{opt-necessary}) is indicated.
Nevertheless, the phase condition~(\ref{phase-ringosc}) yields a
suboptimal solution.

We omit the discussion of solutions of the DAE system~(\ref{ringosc})
reconstructed by the MPDAE solutions using~(\ref{reconstruction}).
Respective numerical results are presented for several test examples
in the previous works~\cite{narayan,pulch05,pulch08a,pulch15}.

\section{Conclusions}
Initial-boundary value problems of MPDAEs were solved numerically using
a method of lines.
The resulting system of DAEs includes an additional constraint
either from an optimisation or from a phase condition.
In the case of semi-explicit DAEs of index one as circuit model, 
we showed that the differential index of the DAEs increases in several cases
determined by the inclusion of either differential variables or
algebraic variables in the additional constraint.
The necessary condition for an optimal solution
including differential variables only is the unique option to
keep the index equal to one in the method of lines.

\clearpage


\end{document}